 \documentclass[12pt,twoside]{article}

 \usepackage{latexsym}
 \usepackage{mathrsfs}
 \usepackage{amsfonts}
 \usepackage{amssymb}
 \usepackage{epsfig}       
 \usepackage{subfigure}    
 \usepackage{graphicx}
 \usepackage{float}

 \def \beginproof{\par\noindent {\bf Proof.}\ }
 \def \endproof{\hbox{ }\hfill$\Box$}

 \def\reff#1{{\rm(\ref{#1})}}

\begin{document}

 \newtheorem{property}{Property}[section]
 \newtheorem{proposition}{Proposition}[section]
 \newtheorem{append}{Appendix}[section]
 \newtheorem{definition}{Definition}[section]
 \newtheorem{lemma}{Lemma}[section]
 \newtheorem{corollary}{Corollary}[section]
 \newtheorem{theorem}{Theorem}[section]
 \newtheorem{algorithm}{Algorithm}[section]
 \newtheorem{remark}{Remark}[section]
 \newtheorem{problem}{Problem}[section]

 \title{A Splitting Augmented Lagrangian Method  for Low Multilinear-Rank Tensor Recovery \thanks{This work was partially supported by National Nature Science Foundation of China (No. 11171252).}}

 \author{Lei Yang\thanks{Department of Mathematics, School of Science, Tianjin University, Tianjin 300072, P.R. China. Email: ylei@tju.edu.cn}\quad
 Zheng-Hai Huang\thanks{Corresponding author. Department of
 Mathematics, School of Science, Tianjin University, Tianjin 300072,
 P.R. China. He is also with the Center for Applied Mathematics of Tianjin
 University. Email: huangzhenghai@tju.edu.cn}\quad
 Yufan Li\thanks{Department of Mathematics, School of Science, Tianjin University, Tianjin 300072, P.R. China. Email: liyufan@tju.edu.cn}}

\date{}

\maketitle

\begin{abstract}
\noindent This paper studies a recovery task of finding a low multilinear-rank tensor that fulfills some linear constraints in the
general settings, which has many applications in computer vision and graphics. This problem is named as the low multilinear-rank tensor recovery problem. The variable splitting technique and convex relaxation technique are used to transform this
problem into a tractable constrained optimization problem. Considering the favorable structure of the problem, we develop a splitting augmented Lagrangian method to solve the resulting problem. The proposed algorithm is easily implemented and its convergence can be proved under some conditions. Some preliminary numerical results on randomly generated and real completion problems show that the proposed algorithm is very effective and robust for tackling the low multilinear-rank tensor completion problem.
 \vspace{5mm}

\noindent {\bf Keywords:}\hspace{2mm} Multilinear-rank; low multilinear-rank tensor recovery; tensor completion;
augmented Lagrangian method
\vspace{2mm}

\noindent {\bf AMS subject classifications:}\hspace{2mm} 90C25, 93C41, 65K05

\end{abstract}

\graphicspath{{images/}}

 \section{Introduction}

 Tensors (or multidimensional arrays) emerge as the higher-order generalization of vectors and matrices.
 More formally, for a positive integer $N$, an $N$-way or $N$th-order real tensor can be regarded as an element of $\mathbb{R}^{n_1 \times \ldots \times n_N}$. Thus, a first-order tensor is a vector and a second-order tensor is a matrix.
 Tensors are applicable in many fields, which involve the multi-way data, such as psychometrics, chemometrics,
 signal processing, computer vision, data mining and elsewhere (see an excellent survey by Kolda and Bader \cite{tb2008}). Tensor rank as
 an important intrinsic characterization of a tensor has been widely discussed based on various tensor decompositions in the literature.
 Considering the rank as a `sparsity' measure, we can recover a tensor by assuming that it is `sparse' and
 solving the following linear constrained problem:
 \begin{eqnarray}\label{problem00}
 \min\limits_{\mathcal{X}}~\Phi(\mathcal{X})  \quad \mbox{\rm s.t.}~
 \mathscr{A}(\mathcal{X})=\mathbf{b}
 \end{eqnarray}
 where $\mathcal{X} \in \mathbb{R}^{n_1 \times \cdots \times n_N}$ is the decision variable, $\Phi(\cdot)$ denotes a kind of tensor rank
 function, $\mathscr{A}$ : $\mathbb{R}^{n_1 \times \ldots \times n_N} \rightarrow \mathbb{R}^p$ with $p\leq\prod^{N}_{i=1}n_i$ is a linear map, and $\mathbf{b} \in \mathbb{R}^p$.
 One of its special cases is the tensor completion problem:
 \begin{eqnarray*}
 \min \limits_{\mathcal{X}}~\Phi(\mathcal{X}) \quad \mbox{\rm s.t.}~ \mathcal{X}_{\Omega} =
 \mathcal{M}_{\Omega},
 \end{eqnarray*}
 where $\mathcal{X}$, $\mathcal{M}$ are $N$-way tensors with identical size in each dimension, and the entries of
 $\mathcal{M}$ in the index set $\Omega$ are given while the remaining entries are missing. It is a missing value estimation
 problem in which a subset of the entries of the tensor $\mathcal{X}$ is given and the unknown entries are to be deduced
 under the low rank assumption. This problem has many applications in computer vision and graphics, e.g., image inpainting \cite{bmgv2000,ng2006}, video inpainting \cite{kr2007}, etc.

 If $\mathcal{X}$ is a second-order tensor, then problem \reff{problem00} reduces to the matrix rank minimization problem \cite{mhs2001,bmp2010}:
 \begin{eqnarray}\label{problem03}
 \min \limits_{X}~\mathrm{rank}(X) \quad\mbox{\rm
 s.t.}\quad\mathscr{A}(X)=\mathbf{b},
 \end{eqnarray}
 where $X \in \mathbb{R}^{m \times n}$ is the decision variable, and the linear map $\mathscr{A}:
 \mathbb{R}^{m \times n} \rightarrow \mathbb{R}^p$ and vector $\mathbf{b} \in \mathbb{R}^p$ are given.
 This problem finds a solution of the lowest rank that fulfills some linear constraints and has a wide range
 of applications in system identification \cite{zl2009}, optimal control \cite{mhs2001} and low-dimensional embedding in Euclidean space \cite{ney1995}, etc. But this problem is NP-hard \cite{j1990}. Recent studies in matrix rank minimization have shown that under certain conditions, one can obtain a solution to the original problem \reff{problem03} via solving a convex relaxation of
 it \cite{bmp2010,eb2009,bwb2011,b2011,ydk2012}. It is well known that the best convex approximation of the rank function over the unit ball of
 matrices with norm less than one is the nuclear norm \cite{bmp2007}, which denotes the sum of the nonzero singular values of matrices. Thus, the closest convex relaxation of problem \reff{problem03} is as follows:
 \begin{eqnarray}\label{problem04}
 \min \limits_{X}~ ||X||_\ast\quad\mbox{\rm
 s.t.}\quad \mathscr{A}(X)=\mathbf{b}.
 \end{eqnarray}
 Several effective algorithms have been proposed for solving \reff{problem04}, such as Fixed Point Continuation with Approximate Singular Value Decomposition (FPCA) by Ma \textit{et al}. \cite{sdl2009}, The Singular Value Thresholding Algorithm (SVT) by Cai \textit{et al}. \cite{jez2008}, The Accelerated Proximal Gradient algorithm (APG) by Toh \textit{et al}. \cite{ks2009}, etc.

 However, for higher-order tensor rank minimization, there are only a few studies due to the complexity of the multi-way data analysis. Liu \textit{et al}. \cite{jppj2009} initialized the study on the tensor completion problem and laid the theoretical foundation of low rank tensor completion. After their pioneering work, convex optimization begins to
 be used on the tensor completion problem. Gandy \textit{et al}. \cite{sbi2011} used the $n$-rank of a tensor as a `sparsity' measure and considered the low-$n$-rank tensor recovery problem, which was more general than the one in \cite{jppj2009}.
 They introduced a tractable convex relaxation of the $n$-rank and proposed two algorithms to solve the low-$n$-rank tensor recovery problem numerically. Recently, Signoretto \textit{et al.} \cite{mqlj2011} showed how convex optimization can be leveraged to deal with a broader set of problems under a low multilinear-rank assumption. The authors proposed a scalable template algorithm based on the Douglas-Rachford splitting technique and explicitly derived the convergence guarantees for proposed algorithm. Encouraging numerical results were reported in \cite{mqlj2011}. More recently, Yang \textit{et al.} \cite{yhs2012} proposed a fixed point iterative method for low $n$-rank tensor pursuit and proved the convergence of the method under some assumptions. The numerical results demonstrated the efficiency of the proposed method, especially for ``easy" problems (with high sampling ratio and low $n$-rank). Some related work can also be found in Marco Signoretto \textit{et al.} \cite{mlj2010,mrbj2011}, and Ryota Tomioka \textit{et al.} \cite{rkh2011,rtkh2011}. Additionally, Zhang \textit{et al}. discussed the exact recovery conditions for tensor $n$-rank minimization via its convex relaxation
 in \cite{mz2012}. Moveover, a third-order tensor recovery problem based on a new tensor rank proposed by Kilmer \textit{et al.} \cite{mknr2011} was also investigated by Yang \textit{et al.} \cite{lzsj2013}.

 In this paper, we use the multilinear-rank of a tensor as a `sparsity' measure and consider the tensor
 multilinear-rank minimization problem. By performing variable splitting, the classical augmented Lagrangian
 method (ALM) can be used to solve the reformulation of the original constrained convex programming directly.
 Moreover, motivated by the work in \cite{ty2011,hty2012,lcm2009}, we take the advantage of the favorable structure
 and develop a splitting augmented Lagrangian method (SALM) as an improvement of the classical ALM.
 The convergence of SALM is also proved under some conditions. We also apply the proposed SALM to solve
 the low multilinear-rank tensor completion problems, denoted by SALM-LRTC.
 Some preliminary numerical results show the fast and accurate recoverability of SALM-LRTC.

 The rest of our paper is organized as follows. In Section 2, we briefly introduce some essential
 notations. Section 3 presents the tensor multilinear-rank minimization problem and its relaxation model.
 In Section 4, the augmented-Lagrangian-type methods, including ALM and SALM, are developed to
 solve the relaxed constrained convex optimization problem. In Section 5, we give the convergence analysis of SALM.
 In Section 6, we report the results of some numerical tests and comparisons among different algorithms.
 Some final remarks are given in the last section.

\section{Notation and Preliminaries}

 In this section, we briefly introduce some essential nomenclatures and notations used in this paper.
 Throughout this paper, scalars are denoted by lowercase letters, e.g., $a,b,c,\cdots$; vectors by
 bold lowercase letters, e.g., $\textbf{a},\textbf{b},\textbf{c}, \cdots$; and matrices by uppercase
 letters, e.g., $A,B,C,\cdots$. An \textit{N}-way tensor is denoted as $\mathcal{X}\in\mathbb{R}^{n_1 \times \ldots \times n_N}$, whose elements are denoted as $x_{i_1 \cdots i_k \cdots i_N}$ where $1\leqslant i_k \leqslant n_k$
 and $1\leqslant k\leqslant N$. Let ${\cal T}:=\mathbb{R}^{n_1 \times \ldots \times n_N}$ denote the set of all the \textit{N}-way
 tensors, and $\mathcal{T}^{N}:=\underbrace{\mathcal{T} \times \mathcal{T} \times \cdots \times \mathcal{T}} \limits_{N}$.

 Next, we present some basic facts about tensors and more details can be found in \cite{tb2008,mqlj2011}. An $N$-way
 tensor $\mathcal{X}$ is rank-$\mathbf{1}$ if it consists of the outer product of $N$ nonzero vectors. Then,
 the linear span of such rank-$\mathbf{1}$ elements forms the vector space ${\cal T}$, which is endowed
 with the inner product
 \begin{eqnarray*}
 \langle\mathcal{X}, \mathcal{Y}\rangle =
 \sum_{i_1=1}^{n_1}\sum_{i_2=1}^{n_2}\cdots\sum_{i_N=1}^{n_N} a_{i_1 i_2 \cdots i_N} b_{i_1 i_2 \cdots i_N}.
 \end{eqnarray*}
 The corresponding (Frobenius-) norm is $\|\mathcal{X}\|_{F}=\sqrt{\langle\mathcal{X},
 \mathcal{X}\rangle}$. When the tensor $\mathcal{X}$ reduces to the
 matrix $X$, $\|X\|_F$ is just the Frobenius norm of the matrix $X$.

 Matricization, also known as unfolding or flattening, is the process of reordering the elements of
 an $N$-way array into a matrix. We use ${X}_{(n)}$ to denote the mode-$n$ unfolding of the tensor
 $\mathcal{X}\in\mathbb{R}^{n_1 \times \ldots \times n_N}$. Specially, the tensor element $(i_1,i_2,\cdots,i_N)$
 is mapped to the matrix element $(i_n,j)$, where
 \begin{eqnarray*}
 j=1+\sum_{k=1,k\neq n}^{N}{(i_k-1)V_k}~~~\mbox{\rm with}~~V_k=\prod_{m=1,m\neq n}^{k-1}{n_m}.
 \end{eqnarray*}
 That is, the ``unfold" operation on a tensor $\mathcal{X}$ is defined as unfold$_n(\mathcal{X}):=
 {X}_{(n)}\in \mathbb{R}^{n_n\times J_n}$ with $J_n=\prod_{k=1,k\neq n}^{N}{n_k}$. The opposite operation
 ``refold" is defined as refold$_{n}({X}_{(n)}):= \mathcal{X}$. The $n$-rank of an $N$-way tensor $\mathcal{X}\in\mathbb{R}^{n_1 \times \ldots \times n_N}$, indicated by rank$(X_{(n)})$, is the rank of ${X}_{(n)}$. A tensor with $r_n = \mathrm{rank}(X_{(n)})$ for $n \in \{ 1, 2, \ldots, N \}$ is called a rank-$(r_1, r_2, \ldots, r_N)$ tensor; the $N$-tuple $(r_1, r_2, \ldots, r_N)$ is called as \textit{multilinear-rank} \footnote{Note that the $n$-rank of $\mathcal{X}$ is defined by the rank of $X_{(n)}$ in \cite{tb2008,mqlj2011} and the $N$-tuple $(r_1, r_2, \ldots, r_N)$ is called as multilinear-rank in \cite{mqlj2011}. While the $N$-tuple $(r_1, r_2, \ldots, r_N)$ is also called as $n$-rank in \cite{sbi2011}. In order to prevent confusion, we adopt the same notation and
 conventions as \cite{tb2008,mqlj2011} in this paper} of $\mathcal{X}$.

 An alternative notion of rank is $\mbox{\rm rank}_{CP}(\mathcal{X})$ \cite{tb2008}, which is defined as
 the smallest number of rank-$\mathbf{1}$ tensors that generate $\mathcal{X}$ as their sum. Whereas for
 second order tensors, $\mathrm{rank}_1(\mathcal{X})=\mathrm{rank}_2(\mathcal{X})=\mathrm{rank}_{CP}(\mathcal{X})$;
 for the general case, it follows that $\mathrm{rank}_n(\mathcal{X}) \leq \mathrm{rank}_{CP}(\mathcal{X})$
 for any $n \in \{ n_1, n_2, \ldots, n_N \}$.

 The $i$-mode (matrix) product of a tensor $\mathcal{X}\in\mathbb{R}^{n_1 \times \cdots \times n_N}$ with a matrix $U \in \mathbb{R}^{L \times n_i}$ is denoted by
 $\mathcal{X} \times_{i} U$ and is of size $n_1 \times \cdots \times n_{i-1} \times L \times n_{i+1} \times \cdots \times n_N$. It can be expressed in terms of unfolded tensors:
 \begin{eqnarray*}
 \mathcal{Y} = \mathcal{X} \times_{i} U \quad \Longleftrightarrow  \quad Y_{(i)}=U X_{(i)}.
 \end{eqnarray*}

 For any matrix $X$, $\|X\|_2$ denotes the operator norm of the matrix $X$; $||X||_\ast=\sum^{n}_{i=1}\sigma_{i}(Y)$
 denotes the nuclear norm of $X$, i.e., the sum of the singular values of $X$. For any vector $x$,
 we use Diag$(x)$ to denote a diagonal matrix with its $i$-th diagonal element being $x_i$.

 Additionally, for $x \in \mathbb{R}_{+}^n$ and $\tau > 0$, the nonnegative vector shrinkage operator
 $s_{\tau}(\cdot)$ is defined as
 \begin{eqnarray*}
 s_{\tau}(x):=\bar{x}\;\; \mbox{\rm with} ~\bar{x}_i
 =\left\{\begin{array}{ll}x_i-\tau,& \mbox{\rm if}\; x_i - \tau > 0, \\
 0,& \mbox{\rm otherwise}. \end{array}\right.
 \end{eqnarray*}
 For $X \in \mathbb{R}^{m \times n}$ and $\tau > 0$, the matrix shrinkage operator $\mathcal{D}_{\tau}(\cdot)$
 is defined as $\mathcal{D}_{\tau}(X):=U\textrm{Diag}(s_{\tau}(\sigma))V^{\top}$, where $X=U\textrm{Diag}(\sigma)V^{\top}$
 is the singular value decomposition of $X$.

 \section{Low Multilinear-Rank Tensor Recovery Problem}

 In this section, we will consider the multilinear-rank tensor recovery problem. Note that the minimization problem
 \reff{problem00} is extended from the matrix (i.e., second-order tensor) case. But unlike matrix, the tensor rank
 is much more complex and non-unique. In fact, $\mbox{\rm rank}_{CP}(\mathcal{X})$ is difficult to handle, as there is no straightforward algorithm to determine rank$_{CP}$ of a specific given tensor. This problem is NP-hard \cite{j1990}. While multilinear-rank is easy to compute, therefore we pay our attention on the multilinear-rank, which is denoted as $\mathrm{rank}_{multi}(\mathcal{X})$ below, in this work. Then, \reff{problem00} becomes the following minimization problem:
 \begin{eqnarray}\label{problem05}
 \min \limits_{\mathcal{X} \in \mathcal{T}} ~\mathrm{rank}_{multi}(\mathcal{X}) \quad \mbox{\rm s.t.}~~\mathscr{A}(\mathcal{X})=\mathbf{b}.
 \end{eqnarray}
 Since the multilinear-rank is the tuple of the ranks of the mode-$n$ unfoldings, \reff{problem05} is actually a
 multiobjective optimization problem. To keep things simple, we will use the sum of the ranks of the different unfoldings
 as the objective function. Thus, we will consider the following \textit{low multilinear-rank tensor recovery problem}:
 \begin{eqnarray}\label{problem2}
 \min\limits_{\mathcal{X} \in \mathcal{T}}~\sum_{i=1}^{N}\mbox{\rm
 rank}({X}_{(i)})\quad \mbox{\rm s.t.}\quad
 \mathscr{A}(\mathcal{X})=\mathbf{b},
 \end{eqnarray}
 where the linear map
 $\mathscr{A}$ : $\mathbb{R}^{n_1 \times \ldots \times n_N} \rightarrow \mathbb{R}^p$ with
 $p\leq\prod^{N}_{i=1}n_i$ and vector $\mathbf{b} \in \mathbb{R}^p$ are given. The corresponding
 \textit{low multilinear-rank tensor completion} is
 \begin{eqnarray}\label{problem3}
 \min\limits_{\mathcal{X} \in \mathcal{T}}~\sum_{i=1}^{N}\mbox{\rm
 rank}({X}_{(i)})\quad \mbox{\rm s.t.}\quad \mathcal{X}_{\Omega} =
 \mathcal{M}_{\Omega},
 \end{eqnarray}
 where $\mathcal{X}$, $\mathcal{M}$ are $N$-way tensors with identical size in each mode, and the entries of
 $\mathcal{M}$ in the index set $\Omega$ are given while the remaining entries are missing.

 The low multilinear-rank tensor recovery problem \reff{problem2} (also its special case \reff{problem3}) is a difficult non-convex problem due to the combination nature of the function rank($\cdot$). Therefore, we also use the nuclear norm as an approximation
 of rank($\cdot$) to get a convex and computationally tractable problem. Then, the nuclear norm
 relaxation of \reff{problem2} is:
 \begin{eqnarray}\label{problem4}
 \min\limits_{\mathcal{X} \in
 \mathcal{T}}~\sum_{i=1}^{N}||{X}_{(i)}||_\ast\quad\mbox{\rm
 s.t.}\quad \mathscr{A}(\mathcal{X})=\mathbf{b}.
 \end{eqnarray}

 Note that problem \reff{problem4} is difficult to solve due to the interdependent nuclear norms. Therefore,
 we perform variable splitting and attribute a separate variable to each unfolding of $\mathcal{X}$.
 Let $\mathcal{Y}_1, \mathcal{Y}_2, \cdots, \mathcal{Y}_N$ be new tensor variables, which are
 equal to the tensor $\mathcal{X}$, i.e., introduce the new variable $ \mathcal{Y}_i \in \mathcal{T}$
 such that $Y_{i, (i)}=X_{(i)}$ for all $i \in \{1, 2, \cdots, N \}$. With these new
 variables $\mathcal{Y}_i$s, we can rephrase \reff{problem4} as follows:
 \begin{eqnarray}\label{problem5}
 &\min \limits_{\mathcal{X}, \mathcal{Y}_i}&~\sum \limits_{i=1}^{N}||Y_{i, (i)}||_\ast   \nonumber \\
 &\mathrm{s.t.}& \mathcal{Y}_i = \mathcal{X}, \quad \forall i \in \{ 1, 2, \cdots, N \}, \\
 &&\mathcal{Y}_i \in \mathcal{T}, ~\mathcal{X} \in \mathcal{B}:= \{\mathcal{X} \in \mathcal{T} | \mathscr{A}(\mathcal{X})=\mathbf{b}\}.  \nonumber
 \end{eqnarray}
 In the following, we design an algorithm to solve \reff{problem5}.

 \section{The Augmented-Lagrangian-Type Methods}

 In this section, the augmented-Lagrangian-type methods are introduced to solve problem \reff{problem5}.

 First, we directly apply the classical ALM (i.e., augmented Lagrangian method) to solve problem \reff{problem5} by noticing that it's a convex constrained optimization problem. As we all know, ALM is one of the state-of-the-art methods for solving constrained optimization problem. It has a pleasing convergence speed with higher accuracy under some rather general conditions. Then, by introducing Lagrange multiplier $\Lambda_i \in\mathcal{T}$ to the equality constraint $\mathcal{Y}_i = \mathcal{X}$ for any $i \in \{ 1,2,\cdots,N\}$, we can give the augmented Lagrangian function  of \reff{problem5}:
 \begin{eqnarray}\label{lagrange}
 \begin{array}{ll}
 \mathcal{L}_A(\mathcal{X}, \mathcal{Y}_1, \cdots, \mathcal{Y}_N, \Lambda_1, \cdots, \Lambda_N, \beta) \\
  \quad:=\sum \limits_{i=1}^{N}||Y_{i, (i)}||_\ast - \sum \limits_{i=1}^N \langle \Lambda_i, \mathcal{X}-\mathcal{Y}_i \rangle
 +\sum \limits_{i=1}^N \frac{\beta}{2} \| \mathcal{X}-\mathcal{Y}_i \|_F^2,
 \end{array}
 \end{eqnarray}
 where $\beta > 0$ is the penalty parameter. Then, the classical ALM applied for
 problem \reff{problem5} is outlined in Algorithm 1.

 \begin{table}[H]
 \centering  \tabcolsep 6pt
 \small{
 \begin{tabular}{l}
 \hline
 \textbf{Algorithm 1} The classical ALM for low multilinear-rank tensor recovery \\
 \hline
 \textbf{Input}: $\mathscr{A},~\mathbf{b},~\Lambda^0,~\beta^0,~\rho \geq 1$  \vspace{0.5mm}\\
 \quad 1: \textbf{while} not converged, \textbf{do}   \vspace{0.5mm}  \\
 \quad 2: \quad Compute $(\mathcal{X}^{k+1}, \mathcal{Y}^{k+1}_1, \cdots, \mathcal{Y}^{k+1}_N)=\mathop{\mathrm{argmin}} \limits_{\mathcal{X} \in \mathcal{B}, \mathcal{Y}_i \in \mathcal{T}} \mathcal{L}_A(\mathcal{X}, \mathcal{Y}_1, \cdots, \mathcal{Y}_N, \Lambda_1^k, \cdots, \Lambda_N^k, \beta^k)$ \vspace{0.5mm} \\
 \quad 3: \quad For $i \in \{1,2,\cdots,N\}$, update $\Lambda^{k+1}_i = \Lambda^k_i-\beta^k(\mathcal{X}^{k+1}-\mathcal{Y}^{k+1}_i)$, \vspace{1mm} \\
 \quad 4: \quad Update $\beta^{k+1}=\rho \beta^k$ \vspace{0.5mm} \\
 \quad 5: \textbf{end while} \vspace{0.5mm} \\
 \textbf{Output}: $(\mathcal{X}^{k}, \mathcal{Y}^{k}_1, \cdots, \mathcal{Y}^{k}_N)$ \\
 \hline
 \end{tabular}}
 \end{table}

 The convergence of Algorithm 1 is easily to derive, one can refer to \cite{m1969,mj1969} for more details. For the implement of the subproblem:
 \begin{eqnarray}\label{subpro}
 (\mathcal{X}^{k+1}, \mathcal{Y}^{k+1}_1, \cdots, \mathcal{Y}^{k+1}_N)
  =\mathop{\mathrm{argmin}} \limits_{\mathcal{X} \in \mathcal{B}, \mathcal{Y}_i \in \mathcal{T}} \mathcal{L}_A(\mathcal{X}, \mathcal{Y}_1, \cdots, \mathcal{Y}_N, \Lambda_1^k, \cdots, \Lambda_N^k, \beta^k)
 \end{eqnarray}
 in Step 2 of Algorithm 1, we may cost much time to solve it to get an optimal solution. Then, inspired by Tao \textit{et al.} \cite{ty2011}, He \textit{et al.} \cite{hty2012} and Lin \textit{et al.} \cite{lcm2009}, taking the favorable structure emerging in both the objective function and constraints of \reff{problem5}
 into consideration, we can split subproblem \reff{subpro}
 into $N+1$ independent minimization subproblems for given $(\Lambda_1^k, \cdots, \Lambda_N^k)$ and $\beta^k$:
\begin{eqnarray}\label{spl-subpro}
 \left\{ \begin{array}{lllrll}
         \mathcal{X}^{k+1} \in \mathop{\mathrm{argmin}} \limits_{\mathcal{X} \in \mathcal{B}} \mathcal{L}_A(\mathcal{X}, \mathcal{Y}^k_1, \mathcal{Y}^k_2, \cdots, \mathcal{Y}^k_N, \Lambda_1^k, \cdots, \Lambda_N^k, \beta^k), \\
         \mathcal{Y}_1^{k+1} \in \mathop{\mathrm{argmin}} \limits_{\mathcal{Y}_1 \in \mathcal{T}} \mathcal{L}_A(\mathcal{X}^{k+1}, \mathcal{Y}_1, \mathcal{Y}^k_2, \cdots, \mathcal{Y}^k_N, \Lambda_1^k, \cdots, \Lambda_N^k, \beta^k), \\
         \mathcal{Y}_2^{k+1} \in \mathop{\mathrm{argmin}} \limits_{\mathcal{Y}_2 \in \mathcal{T}} \mathcal{L}_A(\mathcal{X}^{k+1}, \mathcal{Y}_1^{k+1}, \mathcal{Y}_2, \cdots, \mathcal{Y}^k_N, \Lambda_1^k, \cdots, \Lambda_N^k, \beta^k), \\
         \qquad \quad \vdots \\
         \mathcal{Y}_N^{k+1} \in \mathop{\mathrm{argmin}} \limits_{\mathcal{Y}_N \in \mathcal{T}} \mathcal{L}_A(\mathcal{X}^{k+1}, \mathcal{Y}_1^{k+1}, \mathcal{Y}^{k+1}_2, \cdots, \mathcal{Y}_N, \Lambda_1^k, \cdots, \Lambda_N^k, \beta^k),
         \end{array} \right.
 \end{eqnarray}
 which iterate $(\mathcal{X}^k, \mathcal{Y}^k_1, \cdots, \mathcal{Y}^k_N)$ through solving $N+1$ separable subproblems one by one
 based on the idea of alternative direction method and using the latest information in each subproblem.
 Moreover, it's very encouraging to see that the exact solution of each subproblem in \reff{spl-subpro} can be derived easily.
 In the following, we discuss it separately.

 For the first subproblem in \reff{spl-subpro}, i.e., the minimization of $\mathcal{L}_A$ over $\mathcal{B}$ with respect to the variable $\mathcal{X}$, which is actually a quadratic optimization problem:
 \begin{eqnarray}\label{problem6}
 \min \limits_{\mathcal{X} \in \mathcal{B}} ~\mathcal{L}_A(\mathcal{X})=- \sum \limits_{i=1}^N \langle \Lambda_i, \mathcal{X}-\mathcal{Y}_i \rangle
 +\sum \limits_{i=1}^N \frac{\beta}{2} \| \mathcal{X}-\mathcal{Y}_i \|_F^2.
 \end{eqnarray}
 It is easy to check that \reff{problem6} is equivalent to the following problem:
 \begin{eqnarray}\label{problem7}
 \min \limits_{\mathcal{X} \in \mathcal{B}}~ \left\| \mathcal{X}-\frac{1}{N\beta}\left[\sum \limits_{i=1}^{N} \Lambda_i +
 \sum \limits_{i=1}^{N} \beta \mathcal{Y}_i\right]\right\|_F^2.
 \end{eqnarray}
 By the definition of $\mathcal{B}$ in \reff{problem5}, $\mathcal{B}$ is convex and closed. Then,
 the optimal solution of \reff{problem7} can be achieved by
 \begin{eqnarray*}
 \mathcal{X}^*=\mathcal{P}_{\mathcal{B}}\left( \frac{1}{N\beta}\left[\sum \limits_{i=1}^{N} \Lambda_i +
 \sum \limits_{i=1}^{N} \beta \mathcal{Y}_i \right] \right),
 \end{eqnarray*}
 where $\mathcal{P}_{\mathcal{B}}(\cdot)$ is a projection operator onto set $\mathcal{B}$. Especially,
 for the case of tensor completion problem,
 $\mathcal{B}=\{\mathcal{X} \in \mathcal{T} | \mathcal{X}_{\Omega}=\mathcal{M}_{\Omega}\}$ and we
 can readily calculate the optimal solution $\mathcal{X}^*$ via
 \begin{eqnarray*}
 \mathcal{X}_{i_1 i_2 \cdots i_N}^*=\left\{ \begin{array}{ll}
                       \mathcal{M}_{i_1 i_2 \cdots i_N}, \quad
                       \mathrm{if}~ (i_1, i_2, \cdots, i_N) \in \Omega,  \\
                       \left( \frac{1}{N\beta}\left[\sum \limits_{i=1}^{N} \Lambda_i +
 \sum \limits_{i=1}^{N} \beta \mathcal{Y}_i \right] \right)_{i_1 i_2 \cdots i_N}, \quad \mathrm{otherwise}.
                        \end{array}   \right.
 \end{eqnarray*}

 For the left $N$ subproblems which minimize $\mathcal{L}_A$ over $\mathcal{T}$ with respect to the variable $\mathcal{Y}_i$ for any fixed $i \in \{1, 2, \cdots, N\}$, we fix all variables except $\mathcal{Y}_i$ for given $i$ and then the minimization of $\mathcal{L}_A$ becomes:
 \begin{eqnarray}\label{problem8}
 \begin{array}{cl}
 \min \limits_{\mathcal{Y}_i \in \mathcal{T}} ~\mathcal{L}_A(\mathcal{Y}_i)
 =\|Y_{i,(i)}\|_* - \langle \Lambda_i, \mathcal{X}-\mathcal{Y}_i \rangle
 + \frac{\beta}{2} \| \mathcal{X}-\mathcal{Y}_i \|_F^2.
 \end{array}
 \end{eqnarray}
 Next, we will give an optimal solution to problem \reff{problem8}. Before this, we need the following lemma.
 \begin{lemma}\label{add-lem}
 Let $h_i(\mathcal{Z})=\|Z_{(i)}\|_*$ for any $\mathcal{Z} \in \mathbb{R}^{n_1 \times \cdots \times n_N}$
 and $i \in \{ 1, 2, \cdots, N \}$; and $g(Z)=\|Z\|_*$ for any $Z \in \mathbb{R}^{m \times n}$.
 Then, for any $i \in \{ 1, 2, \cdots, N \}$, the subdifferential of $h_i(\cdot)$ at $\mathcal{Z}$ is the mode-$i$
 refolding of the subdifferential of $\|\cdot\|_*$ at $Z_{(i)}$, i.e.,
 $\partial h_i(\mathcal{Z}) = \mathrm{refold}_{i}(\partial g(Z_{(i)}))$.
 \end{lemma}
 \beginproof
 Firstly, note that $\mathrm{unfold}_i(\cdot)$ and $\mathrm{refold}_i(\cdot)$ are linear one-to-one invertible operators.
 Then, by the definition of subdifferential, for any $i \in \{1, 2, \cdots, N\}$ and any $\mathcal{Q} \in \partial h_i(\mathcal{Z}) \subseteq \mathbb{R}^{n_1 \times \cdots \times n_N}$,
 we have
 \begin{eqnarray}\label{sub01}
 h_i(\mathcal{Z}^{'}) \geq h_i(\mathcal{Z}) + \langle \mathcal{Z}^{'}-\mathcal{Z}, ~\mathcal{Q} \rangle, ~ \forall
 \mathcal{Z}^{'} \in \mathbb{R}^{n_1 \times \cdots \times n_N}.
 \end{eqnarray}
 It is easy to see that \reff{sub01} is equivalent to
 \begin{eqnarray*}
 \|Z^{'}_{(i)}\|_* \geq \|Z_{(i)}\|_* + \langle Z_{(i)}^{'}- Z_{(i)}, ~Q_{(i)} \rangle, ~ \forall
 Z_{(i)}^{'} \in \mathbb{R}^{n_i \times J_i},
 \end{eqnarray*}
 where $J_i = \prod_{k=1,k \neq i}^{N}{n_k}$. Hence, $Q_{(i)}$ is the subdifferential of $\|\cdot\|_*$ at $Z_{(i)}$
 by using the definition of subdifferential again, i.e., $Q_{(i)} \in \partial g(Z_{(i)})$.

 On the other hand, for any $i \in \{1, 2, \cdots, N\}$ and any $Q_{(i)} \in \partial g(Z_{(i)}) \subseteq
 \mathbb{R}^{n_i \times J_i}$ with $J_i = \prod_{k=1,k \neq i}^{N}{n_k}$, we have
 \begin{eqnarray*}
 \|Z^{'}_{(i)}\|_* \geq \|Z_{(i)}\|_* + \langle Z_{(i)}^{'}- Z_{(i)}, ~Q_{(i)} \rangle, ~ \forall
 Z_{(i)}^{'} \in \mathbb{R}^{n_i \times J_i},
 \end{eqnarray*}
 and then
 \begin{eqnarray*}
 h_i(\mathcal{Z}^{'}) \geq h_i(\mathcal{Z}) + \langle \mathcal{Z}^{'}-\mathcal{Z}, ~\mathcal{Q} \rangle, ~ \forall
 \mathcal{Z}^{'} \in \mathbb{R}^{n_1 \times \cdots \times n_N}.
 \end{eqnarray*}
 Thus, $\mathcal{Q}$ is also a subdifferential of $h_i(\cdot)$ at $\mathcal{Z}$, i.e., $\mathcal{Q} \in \partial h_i(\mathcal{Z})$.

 Consequently, we get that $\partial h_i(\mathcal{Z}) = \mathrm{refold}_{i}(\partial g(Z_{(i)}))$ for any
 $i \in \{ 1, 2, \cdots, N \}$.
 \endproof

 \begin{theorem}
 For any given $i \in \{ 1, 2, \cdots, N \}$, an optimal solution to problem \reff{problem8} can be given by
 \begin{eqnarray*}
 \mathcal{Y}_i^*=\mathrm{refold}_{i}\left(\mathcal{D}_{\frac{1}{\beta}}\left(X_{(i)}-\frac{1}{\beta} \Lambda_{i,(i)}\right)\right).
 \end{eqnarray*}
 \end{theorem}
 \beginproof
 By using Lemma \ref{add-lem}, we can easily obtain the result of this theorem in a similar way as \cite[Theorem 4.1]{yhs2012}. So, we omit it here.
 \endproof

 It's worth noticing that the derivation of optimal solution $\mathcal{Y}_i^*$ above contains one application of the matrix shrinkage operator followed by a refolding of resulting matrix into a tensor and the calculation is independent of the choice of $i$.
 Thus, for given $(\Lambda_1^{k}, \cdots, \Lambda_N^{k})$ and $\beta^k$, \reff{spl-subpro} can be easily written as the following more specific form:
 \begin{eqnarray}\label{inspl-iter}
 \left\{ \begin{array}{lll}
         \mathcal{X}^{k+1} :=\mathcal{P}_{\mathcal{B}}\left( \frac{1}{N\beta^k}\left[\sum \limits_{i=1}^{N} \Lambda^k_i + \sum \limits_{i=1}^{N} \beta^k \mathcal{Y}^k_i \right] \right), \\
         \mathcal{Y}_i^{k+1}:=\mathrm{refold}\left(\mathcal{D}_{\frac{1}{\beta^k}}\left(X^{k+1}_{(i)}-\frac{1}{\beta^k} \Lambda^k_{i,(i)}\right)\right),
         ~\mathrm{for}~ i \in \{ 1,2,\cdots,N\}.
         \end{array} \right.
 \end{eqnarray}

 Now, we are ready to present our SALM (i.e., splitting augmented Lagrangian method) for solving \reff{problem5}.

 \begin{table}[H]
 \centering  \tabcolsep 8pt
 \begin{tabular}{l}
 \hline
 \textbf{Algorithm 2} The SALM for low multilinear-rank tensor recovery  \\
 \hline
 \textbf{Input}: $\mathscr{A},~\mathbf{b},~\Lambda^0,~\beta^0,~\rho \geq 1$  \vspace{0.5mm} \\
 \quad 1: \textbf{while} not converged, \textbf{do}   \vspace{0.5mm}  \\
 \quad 2: \quad Compute $\mathcal{X}^{k+1}=\mathcal{P}_{\mathcal{B}}\left( \frac{1}{N\beta^k}\left[\sum \limits_{i=1}^{N}
          \Lambda^k_i+\sum \limits_{i=1}^{N} \beta^k \mathcal{Y}^k_i \right] \right)$  \vspace{1mm} \\
 \quad 3: \quad For $i \in \{1,2,\cdots,N\}$, compute
          $\mathcal{Y}_i^{k+1}=\mathrm{refold}\left(\mathcal{D}_{\frac{1}{\beta^k}}\left(X^{k+1}_{(i)}-
          \frac{1}{\beta^k}\Lambda^{k}_{i,(i)}\right)\right)$  \vspace{1mm}    \\
 \quad 4: \quad For $i \in \{1,2,\cdots,N\}$, update $\Lambda^{k+1}_i = \Lambda^k_i-\beta^k(\mathcal{X}^{k+1}-\mathcal{Y}^{k+1}_i)$  \vspace{1mm} \\
 \quad 5: \quad Update $\beta^{k+1}=\rho \beta^k$   \vspace{0.5mm} \\
 \quad 6: \textbf{end while} \vspace{0.5mm}  \\
 \textbf{Output}: $(\mathcal{X}^{k}, \mathcal{Y}^{k}_1, \cdots, \mathcal{Y}^{k}_N)$ \\
 \hline
 \end{tabular}
 \end{table}

 In fact, Algorithm 2 can be viewed as a more general alternating direction method with three or
 more separable alternative parts, whose convergence is still open in the literature when there are no more conditions given.
 Nevertheless, due to the well structure of problem \reff{problem5}, when some restrictive conditions on $\{ \beta^k \}$
 are assumed, the validity and optimality of Algorithm 2 can be guaranteed, which is discussed in the next section.


 \section{Convergence Results}

 In this section, we discuss the convergence of SALM given in Algorithm 2.
 In order to illustrate our convergence results conveniently, we begin by equivalently restating \reff{problem5} as
 the following constrained problem:
 \begin{eqnarray}\label{equi-pro}
 \begin{array}{ll}
 ~~~~~~~~ \min \limits_{\mathcal{X}, \mathcal{Y}_i \in \mathcal{T}}~f(\mathcal{X})+
 \sum \limits_{i=1}^{N}||Y_{i, (i)}||_\ast   \vspace{1mm}\\
 ~~~~~~~~~~~  \mathrm{s.t.} ~~\mathcal{Y}_i = \mathcal{X}, \quad \forall i \in \{ 1, 2, \cdots, N \},
 \end{array}
 \end{eqnarray}
 where $f:\mathcal{T}\rightarrow\mathbb{R}$ is defined by
 \begin{eqnarray*}
 f(\mathcal{X})= \left\{ \begin{array}{ll}
                                 0, \qquad  ~\mathcal{X} \in \mathcal{B}:= \{\mathcal{X} \in \mathcal{T} | \mathscr{A}(\mathcal{X})=\mathbf{b}\}, \\
                                 +\infty, \quad \mathrm{otherwise}.
                                 \end{array} \right.
 \end{eqnarray*}
 Note that $f$ is an indicator function of $\mathcal{B}$ and it is also a convex
 function on $\mathcal{T}$ since $\mathcal{B}$ is a convex set. Moreover, the subdifferential
 of $f$ at $\mathcal{X}$ is the normal cone of $\mathcal{B}$ at $\mathcal{X}$ \cite{r1970},
 denoted by $N_{\mathcal{B}}(\mathcal{X})$:
 \begin{eqnarray*}
 N_{\mathcal{B}}(\mathcal{X}):= \left\{ \begin{array}{ll}
                                 \{ \mathcal{Z} \in \mathcal{T}: \langle \mathcal{Z}, \mathcal{X}^{'}
                                 -\mathcal{X}\rangle \leq 0, ~\forall \mathcal{X}^{'} \in \mathcal{B} \},
                                  ~\mathcal{X} \in \mathcal{B}, \vspace{0.5mm}  \\
                                \emptyset, \quad
                                 \mathrm{otherwise}.
                                 \end{array} \right.
 \end{eqnarray*}
 Assume that $(\mathcal{X}^*, \mathcal{Y}^*_1, \cdots, \mathcal{Y}^*_N) \in \mathcal{T}^{N+1}$ is an optimal solution of \reff{equi-pro}, then by the optimality condition of \reff{equi-pro}, there must exist Lagrange multipliers $(\Lambda_1^*, \Lambda_2^*, \cdots, \Lambda^*_N) \in \mathcal{T}^{N}$ satisfying the following inclusions:
 \begin{eqnarray}\label{VI}
 \left\{ \begin{array}{lll}
 \sum_{i=1}^N \Lambda_i^* \in N_{\mathcal{B}}(\mathcal{X}^*), ~\mathcal{X}^* \in \mathcal{B}, \vspace{1mm} \\
 0 \in \partial \left(\|Y^*_{i, (i)}\|_* \right) + \Lambda^*_i, ~\forall i \in \{ 1, 2, \cdots, N \}, \vspace{1mm} \\
 \mathcal{X}^*= \mathcal{Y}_i^*, ~\forall i \in \{ 1, 2, \cdots, N \}.
 \end{array} \right.
 \end{eqnarray}
 Therefore, a solution
 of \reff{VI} also yields an optimal solution of \reff{equi-pro} (Hence \reff{problem5}). Throughout this paper,
 we assume that the optimal solution set of \reff{VI} is nonempty.

 In the following, we show some crucial lemmas before giving the main convergence results of SALM.
 For convenience, we denote
 \begin{eqnarray*}
 \mathscr{V}=\left( \begin{array}{cccccc}
 \mathcal{Y}_1 \\ \vdots \\ \mathcal{Y}_N \\ \Lambda_1 \\ \vdots \\ \Lambda_N \end{array} \right), \quad
 \mathscr{V}_{[1]}=\left( \begin{array}{cccccc}
 Y_{1,(1)} \\ \vdots \\  Y_{N,(1)} \\ \Lambda_{1,(1)} \\ \vdots \\ \Lambda_{N,(1)} \end{array} \right) \in \mathbb{R}^{2Nn_1
 \times J_1},
 \end{eqnarray*}
 where $J_1=\prod_{k=2}^{N}{n_k}$. Note that $\mathscr{V}_{[1]}$ is only a notation of abstract vector and is not the mode-1 unfolding of $\mathscr{V}$. For any $\mathscr{V}_{[1]}^{'}$, $\mathscr{V}_{[1]}^{''}$,
 we define
 \begin{eqnarray}\label{define-inpro}
 \langle \mathscr{V}_{[1]}^{'}, \mathscr{V}_{[1]}^{''}\rangle=\sum_{i=1}^N \left(\langle Y_{i,(1)}^{'},Y_{i,(1)}^{''}\rangle +\langle \Lambda_{i,(1)}^{'},\Lambda_{i,(1)}^{''}\rangle\right).
 \end{eqnarray}
 For any $\mathscr{V}$, we define
 \begin{eqnarray}\label{define-norm}
 \| \mathscr{V} \|_{M,1}=\| \mathscr{V}_{[1]} \|_{M}:=\sqrt{\langle\mathscr{V}_{[1]}, M \cdot \mathscr{V}_{[1]} \rangle}
 \end{eqnarray}
 where $M \in \mathbb{R}^{2Nn_1 \times 2Nn_1}$ is a positive definite matrix. It is easy to check that $\| \cdot \|_{M,1}$ satisfies the properties
 of a norm, and we call it the $(M,1)$-norm associated with $\mathscr{V}$. Analogously, we can define the $(M,i)$-norm for any other $i \in \{ 2, 3, \cdots, N\}$. Here, we only use the $(M,1)$-norm in the rest of the paper.

 \begin{lemma}\label{lem-01}
 Suppose that the sequences $\{(\mathcal{X}^k,\mathcal{Y}_1^k,\cdots,\mathcal{Y}_N^k)\}$ and $\{(\Lambda^k_1,\cdots, \Lambda^k_N)\}$ are generated by SALM, and $(\mathcal{X}^*, \mathcal{Y}^*_1, \cdots, \mathcal{Y}^*_N)$ is an optimal solution of \reff{equi-pro} with $(\Lambda_1^*, \cdots, \Lambda^*_N)$ being the corresponding Lagrange multipliers. Then, the following inequalities hold:
 \begin{eqnarray*}
 \left\{ \begin{array}{lll}
 \langle \mathcal{X}^{k+1}-\mathcal{X}^{*}, \sum_{i=1}^N \tilde{\Lambda}_i^{k+1} - \sum_{i=1}^N
 \Lambda_i^{*} \rangle  \geq 0, \\
 \langle \mathcal{Y}_{i}^{k+1}-\mathcal{Y}_{i}^{k}, -\Lambda_{i}^{k+1}-(-\Lambda_{i}^{k}) \rangle \geq 0,
 ~\forall i \in \{1, 2, \cdots, N\},\\
 \langle \mathcal{Y}_i^{k+1}-\mathcal{Y}^{*}_i, -\Lambda^{k+1}_i -(-\Lambda_i^{*})\rangle \geq 0,
 ~\forall i \in \{1, 2, \cdots, N\},
 \end{array} \right.
 \end{eqnarray*}
 where $\tilde{\Lambda}_i^{k+1} = \Lambda^k_i - \beta^k(\mathcal{X}^{k+1}-\mathcal{Y}_i^k)$ for any
 $i \in \{ 1, 2, \cdots, N \}$.
 \end{lemma}
 \beginproof
 Based on the optimality condition of \reff{spl-subpro}, the $k$-th iteration can be characterized by
 the following system:
 \begin{eqnarray*}
 \left\{ \begin{array}{ll}
 \sum \limits_{i=1}^N \Lambda_i^k - \beta^k \sum \limits_{i=1}^N (\mathcal{X}^{k+1}-\mathcal{Y}_i^k) \in N_{\mathcal{B}}(\mathcal{X}^{k+1}), ~\mathcal{X}^{k+1} \in \mathcal{B}, \\
 0 \in \partial \left(\|Y^{k+1}_{i, (i)}\|_* \right) + \Lambda^k_i - \beta^k (\mathcal{X}^{k+1}-\mathcal{Y}_i^{k+1}),~\forall i \in \{ 1, 2, \cdots, N \}, \\
 \Lambda^{k+1}_i = \Lambda^k_i-\beta^k(\mathcal{X}^{k+1}-\mathcal{Y}^{k+1}_i), ~\forall i \in \{1,2,\cdots,N\},
 \end{array} \right.
 \end{eqnarray*}
 which is equivalent to
 \begin{eqnarray}\label{scheme-1}
 \left\{ \begin{array}{ll}
 \sum \limits_{i=1}^N \tilde{\Lambda}_i^{k+1}  \in N_{\mathcal{B}}(\mathcal{X}^{k+1}), ~\mathcal{X}^{k+1} \in \mathcal{B},\\
 -\Lambda^{k+1}_i \in \partial \left(\|Y^{k+1}_{i, (i)}\|_* \right), ~\forall i \in \{ 1, 2, \cdots, N \}.
 \end{array} \right.
 \end{eqnarray}
 Moreover, by the optimality condition \reff{VI}, we have that
 \begin{eqnarray}\label{scheme-2}
 \left\{ \begin{array}{ll}
 \sum \limits_{i=1}^N \Lambda_i^* \in N_{\mathcal{B}}(\mathcal{X}^*), ~\mathcal{X}^* \in \mathcal{B}, \\
 -\Lambda^*_i \in \partial \left(\|Y^*_{i, (i)}\|_*\right), ~\forall i \in \{ 1, 2, \cdots, N \}.  \\
 \end{array} \right.
 \end{eqnarray}

 From the first formulas of \reff{scheme-1} and \reff{scheme-2}, together with the definition of $N_{\mathcal{B}}(\cdot)$, we can obtain that
 \begin{eqnarray*}
 \begin{array}{l}
 ~~~~\langle \mathcal{X}^{k+1}-\mathcal{X}^{*}, \sum \limits_{i=1}^N \tilde{\Lambda}_i^{k+1} - \sum \limits_{i=1}^N
 \Lambda_i^{*} \rangle
 = \langle -\sum \limits_{i=1}^N \tilde{\Lambda}_i^{k+1}, \mathcal{X}^{*}-\mathcal{X}^{k+1} \rangle
 + \langle -\sum \limits_{i=1}^N \Lambda_i^{*}, \mathcal{X}^{k+1}-\mathcal{X}^{*} \rangle \\
 ~~~~~~~~~~~~~~~~~~~~~~~~~~~~~~~~~~~~~~~~~~~\,\geq 0.
 \end{array}
 \end{eqnarray*}

 Note that $\|Y_{i, (i)}\|_*$ can be viewed as a convex function of $\mathcal{Y}_{i}$, by combining the fact that the subdifferential operator of a convex function is monotone \cite{r1970} and the second formulas of \reff{scheme-1} and \reff{scheme-2}, we can get that
 \begin{eqnarray*}
 \langle \mathcal{Y}_{i}^{k+1}-\mathcal{Y}_{i}^{k}, -\Lambda_{i}^{k+1}-(-\Lambda_{i}^{k}) \rangle \geq 0,& \\
 \langle \mathcal{Y}_i^{k+1}-\mathcal{Y}^{*}_i, -\Lambda^{k+1}_i -(-\Lambda_i^{*})\rangle \geq 0,&
 \end{eqnarray*}
 for any $i \in \{ 1, 2, \cdots, N \}$. The proof is complete.
 \endproof

 According to Lemma \textbf{\ref{lem-01}}, we can give the following lemma.
 \begin{lemma}\label{lem-2}
 Suppose that the sequences $\{(\mathcal{X}^k,\mathcal{Y}_1^k,\cdots,\mathcal{Y}_N^k)\}$ and $\{(\Lambda^k_1,\cdots, \Lambda^k_N)\}$ are generated by SALM, $(\mathcal{X}^*, \mathcal{Y}^*_1, \cdots, \mathcal{Y}^*_N)$ is an optimal solution of \reff{equi-pro} with $(\Lambda_1^*, \cdots, \Lambda^*_N)$ being the corresponding Lagrange multipliers, and $\mathscr{V}$, $\mathscr{V}_{[1]}$ are defined as before. Then,
 \begin{itemize}
 \item[(i)] $\langle \mathscr{V}_{[1]}^{k+1}-\mathscr{V}_{[1]}^k, M_k\cdot(\mathscr{V}_{[1]}^{k+1}-\mathscr{V}_{[1]}^*) \rangle \leq 0$.
 \end{itemize}
 Moreover, if $\{\beta^k\}$ is nondecreasing, we have
 \begin{itemize}
 \item[(ii)] \small{$\| \mathscr{V}^{k+1}-\mathscr{V}^{*} \|_{M_{k+1},1}^2 \leq \| \mathscr{V}^{k}-\mathscr{V}^{*} \|_{M_{k},1}^2- \| \mathscr{V}^{k+1}-\mathscr{V}^{k} \|_{M_{k},1}^2$},
 \item[(iii)] $-\sum \limits_{k=0}^{+\infty} \langle \mathscr{V}_{[1]}^{k+1}-\mathscr{V}_{[1]}^k,
           M_k\cdot(\mathscr{V}_{[1]}^{k+1}-\mathscr{V}_{[1]}^*) \rangle < +\infty$,
 \end{itemize}
 where $M_k=\left[\begin{array}{cc}
 I_{Nn_1} &    \\
 & \frac{1}{(\beta^k)^2} I_{Nn_1}
 \end{array}\right] \in \mathbb{R}^{2N n_1 \times~ 2N n_1}$ is a block diagonal matrix with $I_{Nn_1}\in\mathbb{R}^{Nn_1 \times Nn_1}$ being the identity matrix.
 \end{lemma}
 \beginproof
 (i) From the optimality condition, we have $\mathcal{Y}^*_i = \mathcal{X}^*$ for any $i \in \{1,\cdots,N\}$.
 Combining $\Lambda^{k+1}_i - \Lambda_i^{k}=-\beta^k (\mathcal{X}^{k+1}-\mathcal{Y}^{k+1}_i)$ and
 $\Lambda^{k+1}_i - \tilde{\Lambda}_i^{k+1} = \beta^k(\mathcal{Y}^{k+1}_i-\mathcal{Y}^{k}_i)$, we get
 \begin{eqnarray*}
 &&\langle \Lambda^{k+1}_i - \Lambda_i^{k}, \Lambda^{k+1}_i - \Lambda_i^{*} \rangle \\
 &=& -\beta^k \langle \mathcal{X}^{k+1}-\mathcal{Y}^{k+1}_i, \Lambda^{k+1}_i - \Lambda_i^{*} \rangle \\
 &=& -\beta^k \langle \mathcal{X}^{k+1}-\mathcal{X}^{*}, \Lambda^{k+1}_i - \Lambda_i^{*} \rangle
    +\beta^k \langle \mathcal{Y}_i^{k+1}-\mathcal{Y}^{*}_i, \Lambda^{k+1}_i - \Lambda_i^{*} \rangle \\
 &=&  -\beta^k \langle \mathcal{X}^{k+1}-\mathcal{X}^{*}, \Lambda^{k+1}_i - \tilde{\Lambda}_i^{k+1} \rangle -\beta^k \langle \mathcal{X}^{k+1}-\mathcal{X}^{*}, \tilde{\Lambda}_i^{k+1} - \Lambda_i^{*} \rangle \\
  &&  +\beta^k \langle \mathcal{Y}_i^{k+1}-\mathcal{Y}^{*}_i, \Lambda^{k+1}_i - \Lambda_i^{*} \rangle \\
 &=& -(\beta^k)^2 \langle \mathcal{X}^{k+1}-\mathcal{X}^{*}, \mathcal{Y}^{k+1}_i-\mathcal{Y}^{k}_i \rangle -\beta^k \langle \mathcal{X}^{k+1}-\mathcal{X}^{*}, \tilde{\Lambda}_i^{k+1} - \Lambda_i^{*} \rangle \\
  &&    +\beta^k \langle \mathcal{Y}_i^{k+1}-\mathcal{Y}^{*}_i, \Lambda^{k+1}_i - \Lambda_i^{*} \rangle.
 \end{eqnarray*}
 And thus,
 \begin{eqnarray}\label{inequ-3}
 &&\langle \mathscr{V}_{[1]}^{k+1}-\mathscr{V}_{[1]}^k, M_k\cdot(\mathscr{V}_{[1]}^{k+1}-\mathscr{V}_{[1]}^*) \rangle \nonumber\\
 &=& \sum \limits_{i=1}^N \left( \langle Y_{i,(1)}^{k+1}-Y_{i,(1)}^{k}, Y_{i,(1)}^{k+1}-Y_{i,(1)}^{*} \rangle
 +\frac{1}{(\beta^k)^2} \langle \Lambda_{i,(1)}^{k+1}-\Lambda_{i,(1)}^{k}, \Lambda_{i,(1)}^{k+1}-\Lambda_{i,(1)}^{*} \rangle  \right)  \nonumber\\
 &=&\sum \limits_{i=1}^N \left( \langle \mathcal{Y}_{i}^{k+1}-\mathcal{Y}_{i}^{k}, \mathcal{Y}_{i}^{k+1}-\mathcal{Y}_{i}^{*} \rangle
 +\frac{1}{(\beta^k)^2} \langle \Lambda_{i}^{k+1}-\Lambda_{i}^{k}, \Lambda_{i}^{k+1}-\Lambda_{i}^{*} \rangle \right) \nonumber\\
 &=& \sum \limits_{i=1}^N \left( \langle \mathcal{Y}_{i}^{k+1}-\mathcal{Y}_{i}^{k}, \mathcal{Y}_{i}^{k+1} -\mathcal{Y}_{i}^{*} \rangle -
     \langle \mathcal{X}^{k+1}-\mathcal{X}^{*}, \mathcal{Y}^{k+1}_i-\mathcal{Y}^{k}_i \rangle \right)  \nonumber\\
 &&-\frac{1}{\beta^k} \sum \limits_{i=1}^N \langle \mathcal{X}^{k+1}-\mathcal{X}^{*}, \tilde{\Lambda}_i^{k+1} - \Lambda_i^{*} \rangle +\frac{1}{\beta^k} \sum \limits_{i=1}^N \langle \mathcal{Y}_i^{k+1}-\mathcal{Y}^{*}_i, \Lambda^{k+1}_i - \Lambda_i^{*} \rangle \nonumber\\
 &=& \frac{1}{\beta^k}\sum \limits_{i=1}^N \left(\langle\mathcal{Y}_{i}^{k+1}-\mathcal{Y}_{i}^{k}, \Lambda_{i}^{k+1}-\Lambda_{i}^{k}
     \rangle - \langle \mathcal{X}^{k+1}-\mathcal{X}^{*}, \tilde{\Lambda}_i^{k+1} - \Lambda_i^{*} \rangle \right.   \nonumber\\
 &&  \left. + \langle \mathcal{Y}_i^{k+1}-\mathcal{Y}^{*}_i, \Lambda^{k+1}_i - \Lambda_i^{*}
     \rangle \right)  \\
 &\leq& 0,    \nonumber
 \end{eqnarray}
 where the last inequality follows from Lemma \textbf{\ref{lem-01}}.

 (ii) From the definition of $(M,1)$-norm in \reff{define-norm}, it is easy to verify that
 \begin{eqnarray}\label{equ-1}
 \begin{array}{rl}
 \| \mathscr{V}^{k}-\mathscr{V}^{*} \|_{M_{k},1}^2
 =&\| \mathscr{V}^{k}_{[1]}-\mathscr{V}^{*}_{[1]} \|_{M_{k}}^2\\
 =&\sum \limits_{i=1}^N \left( \| \mathcal{Y}_i^{k}-\mathcal{Y}_i^{*} \|_F^2 + \frac{1}{(\beta^k)^2}\| \Lambda_i^{k} - \Lambda_i^{*} \|_F^2 \right).
 \end{array}
 \end{eqnarray}
 Then, by $\beta^{k+1} \geq \beta^k$ and (i), we can derive that
 \begin{eqnarray}\label{inequ-4}
 &&\| \mathscr{V}^{k+1}-\mathscr{V}^{*} \|_{M_{k+1},1}^2  \nonumber \\
 &\leq& \| \mathscr{V}^{k+1}-\mathscr{V}^{*} \|_{M_{k},1}^2   \nonumber \\
 &=& \| \mathscr{V}^{k}-\mathscr{V}^{*} \|_{M_{k},1}^2 - \| \mathscr{V}^{k+1}-\mathscr{V}^{k} \|_{M_{k},1}^2  + 2 \langle \mathscr{V}_{[1]}^{k+1}-\mathscr{V}_{[1]}^k, M_k\cdot(\mathscr{V}_{[1]}^{k+1}-\mathscr{V}_{[1]}^*) \rangle \\
 &\leq& \| \mathscr{V}^{k}-\mathscr{V}^{*} \|_{M_{k},1}^2 - \| \mathscr{V}^{k+1}-\mathscr{V}^{k} \|_{M_{k},1}^2. \nonumber
 \end{eqnarray}

 (iii) By \reff{inequ-4}, we have that
 \begin{eqnarray*}
 - \langle \mathscr{V}_{[1]}^{k+1}-\mathscr{V}_{[1]}^k, M_k\cdot(\mathscr{V}_{[1]}^{k+1}-\mathscr{V}_{[1]}^*) \rangle
 \leq \frac{1}{2}\| \mathscr{V}^{k}-\mathscr{V}^{*} \|_{M_{k},1}^2 - \frac{1}{2}\| \mathscr{V}^{k+1}-\mathscr{V}^{*} \|_{M_{k+1},1}^2.
 \end{eqnarray*}
 Then, it is easy to get that
 \begin{eqnarray*}
 -\sum \limits_{k=0}^{+\infty} \langle \mathscr{V}_{[1]}^{k+1}-\mathscr{V}_{[1]}^k, M_k\cdot(\mathscr{V}_{[1]}^{k+1}-\mathscr{V}_{[1]}^*) \rangle \leq \frac{1}{2}\| \mathscr{V}^{0}-\mathscr{V}^{*} \|_{M_{0},1}^2 < +\infty.
 \end{eqnarray*}

 The proof is complete.
 \endproof

 From Lemma \textbf{\ref{lem-2}} (ii), we can see that if $\{\beta^k\}$ is nondecreasing, the sequence $\{\| \mathscr{V}^{k}-\mathscr{V}^{*} \|_{M_{k},1}^2\}$ is nonincreasing. Based on these three results in Lemma \textbf{\ref{lem-2}}, we have the following lemma immediately, which paves the way towards the convergence of SALM.

 \begin{lemma}\label{lem-3}
 Let $\{\beta^k\}$ be nondecreasing, and the sequences $\{(\mathcal{X}^k,\mathcal{Y}_1^k,\cdots,\mathcal{Y}_N^k)\}$ and $\{(\Lambda^k_1$, $\cdots$, $\Lambda^k_N)\}$ be generated by SALM. Suppose that $(\mathcal{X}^*, \mathcal{Y}^*_1, \cdots, \mathcal{Y}^*_N)$ is an optimal solution of \reff{equi-pro} with $(\Lambda_1^*, \cdots, \Lambda^*_N)$ being the corresponding Lagrange multipliers, and $\mathscr{V}$, $\mathscr{V}_{[1]}$, $M_k$ are defined as before. Then, for any $i \in \{1, 2, \cdots, N\}$,
 \begin{itemize}
 \item[(i)]$\lim \limits_{k \rightarrow \infty} \|\mathcal{X}^{k+1}-\mathcal{Y}^{k+1}_i\|_F=\lim \limits_{k \rightarrow \infty} \|\mathcal{Y}_i^{k+1}-\mathcal{Y}_i^{k}\|_F=\lim \limits_{k \rightarrow \infty} \frac{1}{\beta^k} \|\Lambda_i^{k+1}-\Lambda_i^{k}\|_F=0$; and

 \item[(ii)] the sequences $\{(\mathcal{X}^k,\mathcal{Y}_1^k,\cdots,\mathcal{Y}_N^k)\}$ and $\{(\Lambda_1^k,\cdots,\Lambda_N^k)\}$ are all bounded.
 \end{itemize}
 \end{lemma}
 \beginproof
 (i) From Lemma \textbf{\ref{lem-2}} (ii), we have that
 \begin{eqnarray*}
 \| \mathscr{V}^{k+1}-\mathscr{V}^{k} \|_{M_{k},1}^2
 \leq  \| \mathscr{V}^{k}-\mathscr{V}^{*} \|_{M_{k},1}^2 -\| \mathscr{V}^{k+1}-\mathscr{V}^{*} \|_{M_{k+1},1}^2.
 \end{eqnarray*}
 Furthermore, we can obtain that
 \begin{eqnarray*}
 \sum \limits_{k=0}^{+\infty} \| \mathscr{V}^{k+1}-\mathscr{V}^{k} \|_{M_{k},1}^2 \leq
 \| \mathscr{V}^{0}-\mathscr{V}^{*} \|_{M_{0},1}^2 < +\infty,
 \end{eqnarray*}
 which immediately implies that
 $\lim \limits_{k \rightarrow \infty} \| \mathscr{V}^{k+1}-\mathscr{V}^{k} \|_{M_{k},1}^2 = 0.$
 Recall that
 \begin{eqnarray*}
 \| \mathscr{V}^{k+1}-\mathscr{V}^{k} \|_{M_{k},1}^2
 =\sum \limits_{i=1}^N \left( \| \mathcal{Y}_i^{k+1}-\mathcal{Y}_i^{k} \|_F^2
  + \frac{1}{(\beta^k)^2}\| \Lambda_i^{k+1} - \Lambda_i^{k} \|_F^2 \right),
 \end{eqnarray*}
 we have that for any $i \in \{1, 2, \cdots, N\}$,
 \begin{eqnarray*}
 \lim \limits_{k \rightarrow \infty} \|\mathcal{Y}_i^{k+1}-\mathcal{Y}_i^{k}\|_F=0, ~~~
 \lim \limits_{k \rightarrow \infty} \frac{1}{\beta^k} \|\Lambda_i^{k+1}-\Lambda_i^{k}\|_F=0,
 \end{eqnarray*}
 which, together with $\mathcal{X}^{k+1}-\mathcal{Y}^{k+1}_i=-\frac{1}{\beta^k}(\Lambda_i^{k+1} - \Lambda_i^{k})$, shows that
 \begin{eqnarray}\label{lim-2}
 \lim \limits_{k \rightarrow \infty}\|\mathcal{X}^{k+1}-\mathcal{Y}^{k+1}_i\|_F=\lim \limits_{k \rightarrow \infty}\frac{1}{\beta^k}\|\Lambda_i^{k+1} - \Lambda_i^{k}\|_F =0.
 \end{eqnarray}

 (ii)  By \reff{equ-1} and Lemma \textbf{\ref{lem-2}} (ii), we can easily obtain that
 \begin{eqnarray*}
 \sum \limits_{i=1}^N \| \mathcal{Y}_i^{k}-\mathcal{Y}_i^{*} \|_F^2 \leq \| \mathscr{V}^{k}-\mathscr{V}^{*} \|_{M_{k},1}^2
     \leq\| \mathscr{V}^{0}-\mathscr{V}^{*} \|_{M_{0},1}^2 < +\infty,
 \end{eqnarray*}
 which shows that $ \| \mathcal{Y}_i^{k}-\mathcal{Y}_i^{*} \|_F < +\infty$ for any $i \in \{1, 2, \cdots, N\}$.
 Then, the sequence $\{ \mathcal{Y}_i^{k} \}$ is bounded. Moreover, the boundedness of $\{\mathcal{Y}_i^{k}\}$ and (i) imply that $\{ \mathcal{X}^{k} \}$ is bounded. Additionally, note that
 \begin{eqnarray*}
 -\Lambda^{k+1}_i \in \partial \left(||Y^{k+1}_{i, (i)}||_\ast \right), \quad i \in \{1, 2, \cdots, N\},
 \end{eqnarray*}
 we have that the sequence $\{ \Lambda^{k}_i \}$ is bounded since the dual norm of $\|\cdot\|_*$ is $\|\cdot\|_2$.
 \endproof

 Now, we are ready to give the convergence of SALM.
 \begin{theorem}\label{theo-1}
 Let $\{\beta^k\}$ be nondecreasing and $\sum_{k=0}^{+\infty} \frac{1}{\beta^k} = +\infty$. Suppose that the sequences $\{(\mathcal{X}^k,\mathcal{Y}_1^k,\cdots,\mathcal{Y}_N^k)\}$ and $\{(\Lambda_1^k,\cdots,\Lambda_N^k)\}$ are generated by SALM, and $\mathscr{V}$, $\mathscr{V}_{[1]}$, $M_k$ are defined as before. Then, the sequence $\{(\mathcal{X}^k, \mathcal{Y}^k_1, \cdots, \mathcal{Y}^k_N)\}$ converges to an optimal solution of \reff{equi-pro} (Hence \reff{problem5}).
 \end{theorem}
 \beginproof
  From Lemma \textbf{\ref{lem-3}} (ii), the sequences $\{(\mathcal{X}^k,\mathcal{Y}_1^k,\cdots,\mathcal{Y}_N^k)\}$ and $\{(\Lambda^k_1,\cdots, \Lambda^k_N)\}$
 are bounded, and thus, exist convergent subsequences. Then, by \reff{lim-2}, for any $i \in \{1,2,\cdots,N\}$,
 \begin{eqnarray*}
 \lim \limits_{k \rightarrow +\infty} \|\mathcal{X}^{k}-\mathcal{Y}^{k}_i\|_F=0,
 \end{eqnarray*}
 which, together with the fact that $\mathcal{B}$ is compact, shows that any accumulation point of sequence $\{(\mathcal{X}^k,\mathcal{Y}_1^k,\cdots,\mathcal{Y}_N^k)\}$ is a feasible solution of \reff{equi-pro}.

 Next, we show that some accumulation point of the sequence $\{(\mathcal{X}^k,\mathcal{Y}_1^k,\cdots,\mathcal{Y}_N^k)\}$ is an optimal solution of \reff{equi-pro}.
 Let $F^*$ denote the optimal objective value of \reff{equi-pro}, and suppose that $(\mathcal{X}^*, \mathcal{Y}^*_1, \cdots, \mathcal{Y}^*_N)$ is an optimal solution of \reff{equi-pro} with $(\Lambda_1^*, \Lambda_2^*, \cdots, \Lambda^*_N)$ being the corresponding Lagrange multipliers. Note that $\|\cdot\|_*$ is convex, which, together with the optimal condition:
 $\sum_{i=1}^N \tilde{\Lambda}_i^{k} \in N_{\mathcal{B}}(\mathcal{X}^{k})$ and $-\Lambda^{k}_i \in \partial (||Y^{k}_{i, (i)}||_\ast)$, $\mathcal{X}^{*}=\mathcal{Y}_i^{*}$ for any $i\in \{1,2,\cdots,N\}$, implies that
 \begin{eqnarray}\label{inequ-5}
 && \sum \limits_{i=1}^{N}||Y^k_{i, (i)}||_\ast    \nonumber\\
 &\leq & \sum \limits_{i=1}^{N} \left( ||Y^*_{i, (i)}||_\ast- \langle -\Lambda^{k}_i, \mathcal{Y}^*_i - \mathcal{Y}^k_i \rangle \right)
         -\langle \sum \limits_{i=1}^N \tilde{\Lambda}_i^{k}, \mathcal{X}^*-\mathcal{X}^k \rangle    \nonumber\\
 &=& F^*+\sum \limits_{i=1}^N \langle -\Lambda_{i}^{*}-(-\Lambda_{i}^{k}), \mathcal{Y}_{i}^{*}-\mathcal{Y}_{i}^{k} \rangle
     + \sum \limits_{i=1}^N \langle \Lambda_i^{*}-\tilde{\Lambda}_i^{k}, \mathcal{X}^{*}-\mathcal{X}^{k} \rangle  \nonumber\\
 & &  +\sum \limits_{i=1}^N \langle \Lambda_{i}^{*}, \mathcal{Y}_{i}^{*}-\mathcal{Y}_{i}^{k} \rangle
    - \sum \limits_{i=1}^N \langle \Lambda_i^{*}, \mathcal{X}^{*}-\mathcal{X}^{k} \rangle    \nonumber\\
 &=& F^*+\sum \limits_{i=1}^N \langle -\Lambda_{i}^{k}-(-\Lambda_{i}^{*}), \mathcal{Y}_{i}^{k}-\mathcal{Y}_{i}^{*} \rangle
     + \sum \limits_{i=1}^N \langle \tilde{\Lambda}_i^{k}-\Lambda_i^{*}, \mathcal{X}^{k}-\mathcal{X}^{*} \rangle \nonumber\\
 & &  + \sum \limits_{i=1}^N \langle \Lambda_{i}^{*}, \mathcal{X}_{i}^{k}-\mathcal{Y}_{i}^{k} \rangle.
 \end{eqnarray}
 By \reff{inequ-3} and Lemma \textbf{\ref{lem-2}} (iii),
 \begin{eqnarray*}
 &&\sum \limits_{k=0}^{+\infty} \frac{1}{\beta^k}\left( \sum \limits_{i=1}^N \langle \mathcal{Y}_i^{k+1}- \mathcal{Y}^{*}_i, -\Lambda^{k+1}_i
   -(-\Lambda_i^{*}) \rangle + \sum \limits_{i=1}^N \langle \mathcal{X}^{k+1}-\mathcal{X}^{*}, \tilde{\Lambda}_i^{k+1} - \Lambda_i^{*} \rangle \right) \\
 &\leq& -\sum \limits_{k=0}^{+\infty} \langle \mathscr{V}_{[1]}^{k+1}-\mathscr{V}_{[1]}^k, M\cdot(\mathscr{V}_{[1]}^{k+1}-\mathscr{V}_{[1]}^*) \rangle < +\infty.
 \end{eqnarray*}
 As $\sum_{k=0}^{+\infty} \frac{1}{\beta^k} = +\infty$, there must exist a subsequence
 $\{(\mathcal{X}^{k_j}, \mathcal{Y}^{k_j}_1, \cdots, \mathcal{Y}^{k_j}_N)\}$ such that
 \begin{eqnarray*}
 \sum \limits_{i=1}^N \left(\langle \mathcal{Y}_i^{k_j}- \mathcal{Y}^{*}_i, -\Lambda^{k_j}_i -(-\Lambda_i^{*}) \rangle +
 \langle \mathcal{X}^{k_j}-\mathcal{X}^{*}, \tilde{\Lambda}_i^{k_j} -\Lambda_i^{*} \rangle \right) \longrightarrow 0
 ~~ \mathrm{as} ~ j\rightarrow +\infty.
 \end{eqnarray*}
 Here, we assume that the subsequence $\{(\mathcal{X}^{k_j}, \mathcal{Y}^{k_j}_1, \cdots, \mathcal{Y}^{k_j}_N)\}$ converges to $(\hat{\mathcal{X}}$, $\hat{\mathcal{Y}_1}$, $\cdots$, $\hat{\mathcal{Y}_N})$. If not, there must exist a convergent subsequence
 of  $\{(\mathcal{X}^{k_j}, \mathcal{Y}^{k_j}_1, \cdots, \mathcal{Y}^{k_j}_N)\}$ since it is bounded by Lemma \textbf{\ref{lem-3}} (ii), and we can denote this subsequence as
 $\{(\mathcal{X}^{k_j}, \mathcal{Y}^{k_j}_1, \cdots, \mathcal{Y}^{k_j}_N)\}$ again without loss of generality.
 By \reff{lim-2}, we have that for any $i \in \{1,2,\cdots,N\}$,
 \begin{eqnarray*}
 \lim \limits_{j \rightarrow +\infty} \|\mathcal{X}^{k_j}-\mathcal{Y}^{k_j}_i\|_F=0 \quad \Longrightarrow \quad \|\hat{\mathcal{X}}-\hat{\mathcal{Y}_i}\|_F=0,
 \end{eqnarray*}
 and furthermore,
 \begin{eqnarray}\label{lim3}
 \hat{\mathcal{X}}=\hat{\mathcal{Y}_i}  \quad \Longrightarrow \quad \lim \limits_{j \rightarrow +\infty} (\mathcal{X}^{k_j}-\mathcal{Y}^{k_j}_i)=0.
 \end{eqnarray}
 Then, for \reff{inequ-5}, we take the limit and get
 \begin{eqnarray*}
 \sum \limits_{i=1}^{N}||\hat{Y}_{i, (i)}||_\ast=\lim \limits_{j \rightarrow +\infty} \sum \limits_{i=1}^{N}||Y^{k_j}_{i, (i)}||_\ast \leq F^*.
 \end{eqnarray*}
 So the limit point $(\hat{\mathcal{X}},\hat{\mathcal{Y}_1},\cdots,\hat{\mathcal{Y}_N})$ of $(\mathcal{X}^{k_j}, \mathcal{Y}^{k_j}_1, \cdots, \mathcal{Y}^{k_j}_N)$ is an optimal solution of \reff{equi-pro}.

 From Lemma \textbf{\ref{lem-3}} (ii), $\{(\Lambda^{k_j}_1,\cdots, \Lambda^{k_j}_N)\}$ is bounded.
 So, we can directly assume that the subsequence $\{(\Lambda^{k_j}_1,\cdots, \Lambda^{k_j}_N)\}$ converges to
 $(\hat{\Lambda}_1,\cdots, \hat{\Lambda}_N)$ without loss of generality.
 Thus, we have that
 \begin{eqnarray*}
 \| \mathscr{V}^{k_j}-\hat{\mathscr{V}} \|_{M_{k_j},1}^2
 =\sum \limits_{i=1}^N \left( \| \mathcal{Y}_i^{k_j}-\hat{\mathcal{Y}_i} \|_F^2
     + \frac{1}{(\beta^{k_j})^2}\| \Lambda_i^{k_j} - \hat{\Lambda_i} \|_F^2 \right) \rightarrow 0
 ~\mathrm{as}~ j\rightarrow +\infty.
 \end{eqnarray*}
 By Lemma \textbf{\ref{lem-2}} (ii), $\{\| \mathscr{V}^{k}-\hat{\mathscr{V}} \|_{M_{k},1}^2\}$ is nonincreasing, then we can derive that
 $$
 \| \mathscr{V}^{k}-\hat{\mathscr{V}} \|_{M_{k},1}^2 \rightarrow 0, \quad k\rightarrow +\infty,
 $$
 which indicates that $\lim \limits_{k \rightarrow +\infty} \mathcal{Y}_i^{k}=\hat{\mathcal{Y}_i}$ for any $i \in \{ 1, 2, \cdots, N \}$.
 By \reff{lim-2} and $\hat{\mathcal{X}}=\hat{\mathcal{Y}_i}$, we see that $\lim \limits_{k \rightarrow +\infty} \mathcal{X}^{k}=\hat{\mathcal{X}}$.

 Consequently, the consequence $\{(\mathcal{X}^k, \mathcal{Y}^k_1, \cdots, \mathcal{Y}^k_N)\}$ converges to $(\hat{\mathcal{X}},\hat{\mathcal{Y}_1},\cdots,\hat{\mathcal{Y}_N})$, which is an optimal solution of \reff{equi-pro}
 (Hence \reff{problem5}). This completes the proof.
 \endproof

 \section{ Numerical Experiments }

 In this section, we apply SALM to solve low multilinear-rank tensor completion problems, which is
 denoted by SALM-LRTC, and evaluate its empirical performance both on simulated and real world data with the missing data.
 We also compare it with the latest tensor completion algorithms, including
 FP-LRTC (fixed point continuation method for low $n$-rank tensor completion) \cite{yhs2012},
 TENSOR-HC (hard completion) \cite{mqlj2011}, ADM-CON (ADMM for the ``Constraint"
 approach) \cite{rkh2011} and ADM-TR(E) (alternative direction method algorithm for low-$n$-rank tensor recovery)
 \cite{sbi2011}. All numerical experiments are run in Matlab 7.11.0 on a HP Z800 workstation with an Intel Xeon(R) 3.33GHz CPU and 48GB of RAM.

 \subsection{Implementation Details}

 \textit{Problem settings}. The random low multilinear-rank tensor completion problems without noise we consider in our numerical experiments are generated as in \cite{sbi2011, rkh2011, yhs2012}.
 For creating a tensor $\mathcal{M}\in\mathbb{R}^{n_1 \times \ldots \times n_N}$ with rank $(r_1,r_2,\cdots,r_N)$, we first generate a core tensor
 $\mathcal{S}\in \mathbb{R}^{r_1\times\cdots\times r_N}$ with i.i.d. Gaussian entries ($\sim\mathcal{N}(0,1)$).
 Then, we generate matrices $U_{1},\cdots,U_{N}$, with $U_{i}\in\mathbb{R}^{n_i\times r_i}$ whose entries
 are i.i.d. from  $\mathcal{N}(0,1)$ and set
 \begin{eqnarray*}
 \mathcal{M}:=\mathcal{S}\times_1U_{1}\times_2\cdots\times_N U_{N}.
 \end{eqnarray*}
 With this construction, the multilinear-rank of $\mathcal{M}$ equals $(r_1,r_2,\cdots,r_N)$ almost surely.

 We also conduct numerical experiments on random low multilinear-rank tensor completion problems with
 noisy data. For the noisy random low multilinear-rank tensor completion problems, the tensor
 $\mathcal{M}\in\mathbb{R}^{n_1 \times \ldots \times n_N}$ is corrupted
 by a noise tensor $\mathcal{E}\in\mathbb{R}^{n_1 \times \ldots \times n_N}$ with independent normally distributed entries.
 Then, $\mathcal{M}$ is taken to be
 \begin{eqnarray}\label{noise}
 \mathcal{M}:=\bar{\mathcal{M}}+\sigma \mathcal{E}=\mathcal{S}\times_1U_{1}\times_2\cdots\times_N U_{N}+\sigma \mathcal{E}.
 \end{eqnarray}

 We use $sr$ to denote the sampling ratio, i.e., a percentage $sr$ of the entries to be known and choose the
 support of the known entries uniformly at random among all supports of size $sr\left(\prod^N_{i=1}n_i\right)$.
 The values and the locations of the known entries of $\mathcal{M}$ are used as input for the algorithms.\\

 \textit{Singular value decomposition}. Computing singular value decomposition (SVD) is the main computational cost of
 solving low multilinear-rank tensor completion problems. So using the Lanczos algorithm \cite{b1992} like PROPACK \cite{propack}
 for computing only partial singular value decomposition is used to speed up the calculation in FP-LRTC, ADM-CON and ADM-TR(E).
 However, to use PROPACK, one have to predict the number of singular values to compute. Although there has been many
 heuristics to choose the predetermined number for matrix or tensor completion \cite{yhs2012, lcm2009, ty2011, rkh2011, sbi2011, ks2009, jez2008, sdl2009}, it is still difficult to give an optimal strategy. In fact, the algorithms may not achieve very good recoverability  under the inappropriate prediction. Hence, for good recoverability and simplicity,
 we use the matlab command $[U,S,V] = \mathrm{svd}(X,'\mathrm{econ}')$ instead of PROPACK to compute full SVD in our
 algorithm. The same command is also used in TENSOR-HC. \\

 \textit{Evaluation criterion and stopping criterion}. For random low multilinear-rank tensor completion problems
 without noise, we report the relative error
 \begin{eqnarray*}
 \mathrm{rel.err}:=\frac{||\mathcal{X}_{\mathrm{sol}}-\mathcal{M}||_F}{||\mathcal{M}||_F}
 \end{eqnarray*}
 to estimate the closeness of $\mathcal{X}_{\mathrm{sol}}$ to $\mathcal{M}$, where $\mathcal{X}_{\mathrm{sol}}$ is the
 ``optimal" solution to \reff{problem4} produced by the algorithms and $\mathcal{M}$ is the original tensor.

 For random low multilinear-rank tensor completion problems with noisy data, we measure the performance based on
 the normalized root mean square error (NRMSE) \cite{mqlj2011} on the complementary set $\Omega^{c}$:
 \begin{eqnarray*}
 \mathrm{NRMSE}(\mathcal{X}^{opt}, \bar{\mathcal{M}}):=\frac{||\mathcal{X}^{opt}_{\Omega^{c}}-\bar{\mathcal{M}}_{\Omega^{c}}||_F}
 {\left(\mathrm{max}(\bar{\mathcal{M}}_{\Omega^{c}})-\mathrm{min}(\bar{\mathcal{M}}_{\Omega^{c}})\right)\sqrt{|\Omega^{c}|}}
 \end{eqnarray*}
 where $\bar{\mathcal{M}}$ is as in \reff{noise} and $|\Omega^{c}|$ denotes the cardinality of $\Omega^{c}$.

 The stopping criterion we use for SALM-LRTC in all our numerical experiments is as follows:
 \begin{eqnarray*}
 \frac{\| \mathcal{X}^{k+1}-\mathcal{X}^k \|_F}{\mathrm{max}\{1, \|\mathcal{X}^k\|_F\}} < \mathrm{Tol},
 \end{eqnarray*}
 where Tol is a moderately small number, since when $\mathcal{X}^k$ gets close to an optimal solution $\mathcal{X}^{opt}$, the distance between $\mathcal{X}^k$ and $\mathcal{X}^{k+1}$ should become very small.   \\

 \textit{Updating $\beta^k$}. In our numerical experiments, to keep things simple, we update $\beta^k$
 adaptively as follows:
 \begin{eqnarray*}
 \beta^{k+1}=\left\{\begin{array}{ll}
 \rho \beta^k,& \mbox{\rm if}~~ \frac{\| \mathcal{X}^{k+1}-\mathcal{X}^k \|_F}{\mathrm{max}\{1, \|\mathcal{X}^k\|_F\}}
 \leq \varepsilon,  \\
 \beta^k,& \mbox{\rm otherwise}, \end{array}\right.
 \end{eqnarray*}
 where $\rho >1$ and $\varepsilon$ is a small positive number. This rule indicates that $\beta^k$ is increased by
 a constant factor $\rho$ when there is a slow change during iterations and it is also used in ADM-TR(E). \\

 \textit{Choice of parameters}. Throughout the experiments, we choose the initial iterate to be $\mathcal{X}^0=0$ and
 set $\beta^0=0.1$, $\rho=5$, $\mathrm{Tol}=10^{-8}$. The parameters in updating $\beta^k$ is set as follows:
 if $sr > 0.5$, we set $\varepsilon=10^{-3}$; otherwise, we set $\varepsilon=10^{-4}$, where $sr$ is the sampling ratio.

 \subsection{Numerical Simulation}

 In this part, we test some randomly created problems to illustrate the recoverability and convergence properties of
 SALM-LRTC and provide the comparisons with other algorithms. All results are average values of 10 independent trials.

 In FP-LRTC, we set $\mu_1=1$, $\tau=10$, $\theta_{\mu}=1-sr$, $\bar{\mu}=1\times10^{-8}$, $\varepsilon = 10^{-2}$.
 In TENSOR-HC, we set the regularization parameters $\lambda_{i}$, $i \in \{1,2,\cdots,N\}$ to $1$ and $\tau$ to $10$.
 The parameter $\lambda$ in ADM-CON is set to 0 and the code of ADM-CON is downloaded from http://www.ibis.t.u-tokyo.ac.jp/RyotaTomioka/Softwares/Tensor.html.
 In ADM-TR(E), the parameters are set to $c_{\beta}=5, c_{\lambda}=5, \beta =1, \lambda = N$.
 The above parameter settings can be seen in \cite{yhs2012, mqlj2011, sbi2011, rkh2011} for more details.

 \begin{table}[H]
 \centering \tabcolsep 6pt
 \begin{tabular}{|llcr|llcr|}
 \multicolumn{8}{c}{\scriptsize{\textbf{Table 1.} Numerical results on low multilinear-rank tensor completion without noise}} \vspace{1mm} \\
 \hline
 \scriptsize{Method} & \scriptsize{iter} & \scriptsize{rel.err} & \scriptsize{time} & \scriptsize{Method} & \scriptsize{iter} & \scriptsize{rel.err} & \scriptsize{time} \\
 \hline
 \multicolumn{8}{|c|}{\scriptsize{$\mathcal{T}=\mathbb{R}^{50\times50\times50}$, $r=(9,9,3)$}}  \\
 \hline
 \multicolumn{4}{|c|}{\scriptsize{$sr=0.3$}} & \multicolumn{4}{c|}{\scriptsize{$sr=0.6$}} \\
 \hline
 \scriptsize{SALM-LRTC}	&\scriptsize{70}&\scriptsize{8.45e-8}&\scriptsize{3.85}& \scriptsize{SALM-LRTC} &\scriptsize{35}&\scriptsize{8.35e-9}&\scriptsize{2.06}	 \vspace{-1.5mm}\\	
 \scriptsize{FP-LRTC}	&\scriptsize{520}&\scriptsize{5.77e-8}&\scriptsize{30.23}& \scriptsize{FP-LRTC} &\scriptsize{105}&\scriptsize{3.09e-8}&\scriptsize{2.03}		 \vspace{-1.5mm}\\
 \scriptsize{TENSOR-HC}	&\scriptsize{76}&\scriptsize{3.81e-8}&\scriptsize{6.75}& \scriptsize{TENSOR-HC} &\scriptsize{59}&\scriptsize{2.24e-8}&\scriptsize{5.19}	 \vspace{-1.5mm}\\	
 \scriptsize{ADM-CON}	&\scriptsize{150}&\scriptsize{3.37e-8}&\scriptsize{13.14}& \scriptsize{ADM-CON} &\scriptsize{64}&\scriptsize{2.69e-8}&\scriptsize{5.52}		 \vspace{-1.5mm}\\
 \scriptsize{ADM-TR(E)}	&\scriptsize{422}&\scriptsize{2.62e-7}&\scriptsize{39.13}& \scriptsize{ADM-TR(E)} &\scriptsize{214}&\scriptsize{1.40e-8}&\scriptsize{22.26}	 \vspace{0mm}\\
 \hline
 \multicolumn{8}{|c|}{\scriptsize{$\mathcal{T}=\mathbb{R}^{100\times100\times50}$, $r=(10,10,5)$}}   \\
 \hline
 \multicolumn{4}{|c|}{\scriptsize{$sr=0.3$}} & \multicolumn{4}{c|}{\scriptsize{$sr=0.6$}} \\
 \hline
 \scriptsize{SALM-LRTC}	&\scriptsize{67}&\scriptsize{9.97e-9}&\scriptsize{22.96}& \scriptsize{SALM-LRTC} &\scriptsize{33}&\scriptsize{4.77e-9}&\scriptsize{11.28}	 \vspace{-1.5mm}\\	
 \scriptsize{FP-LRTC}	&\scriptsize{520}&\scriptsize{2.78e-8}&\scriptsize{39.90}& \scriptsize{FP-LRTC} &\scriptsize{105}&\scriptsize{4.44e-9}&\scriptsize{9.10}		 \vspace{-1.5mm} \\
 \scriptsize{TENSOR-HC}	&\scriptsize{54}&\scriptsize{4.08e-8}&\scriptsize{25.39}& \scriptsize{TENSOR-HC} &\scriptsize{35}&\scriptsize{2.57e-8}&\scriptsize{16.47}	 \vspace{-1.5mm} \\	
 \scriptsize{ADM-CON}	&\scriptsize{149}&\scriptsize{1.85e-8}&\scriptsize{69.29}& \scriptsize{ADM-CON} &\scriptsize{65}&\scriptsize{9.56e-9}&\scriptsize{31.38}		 \vspace{-1.5mm} \\
 \scriptsize{ADM-TR(E)}	&\scriptsize{380}&\scriptsize{1.92e-7}&\scriptsize{175.59}& \scriptsize{ADM-TR(E)} &\scriptsize{186}&\scriptsize{6.86e-9}&\scriptsize{84.57}	 \vspace{0mm} \\
 \hline
 \multicolumn{8}{|c|}{\scriptsize{$\mathcal{T}=\mathbb{R}^{20\times20\times30\times30}$, $r=(4,4,4,4)$}} \\
 \hline
 \multicolumn{4}{|c|}{\scriptsize{$sr=0.3$}} & \multicolumn{4}{c|}{\scriptsize{$sr=0.6$}} \\
 \hline
 \scriptsize{SALM-LRTC}	&\scriptsize{74}&\scriptsize{2.25e-8}&\scriptsize{15.66}& \scriptsize{SALM-LRTC} &\scriptsize{35}&\scriptsize{5.11e-9}&\scriptsize{7.38}	 \vspace{-1.5mm} \\	
 \scriptsize{FP-LRTC}	&\scriptsize{520}&\scriptsize{4.20e-8}&\scriptsize{100.29}& \scriptsize{FP-LRTC} &\scriptsize{210}&\scriptsize{6.96e-9}&\scriptsize{16.53}		 \vspace{-1.5mm} \\
 \scriptsize{TENSOR-HC}	&\scriptsize{52}&\scriptsize{4.49e-7}&\scriptsize{17.06}& \scriptsize{TENSOR-HC} &\scriptsize{36}&\scriptsize{7.82e-8}&\scriptsize{11.92}	 \vspace{-1.5mm} \\	
 \scriptsize{ADM-CON}	&\scriptsize{173}&\scriptsize{1.21e-7}&\scriptsize{53.12}& \scriptsize{ADM-CON} &\scriptsize{72}&\scriptsize{8.17e-8}&\scriptsize{22.12}		 \vspace{-1.5mm} \\
 \scriptsize{ADM-TR(E)}	&\scriptsize{445}&\scriptsize{2.57e-7}&\scriptsize{147.89}& \scriptsize{ADM-TR(E)} &\scriptsize{228}&\scriptsize{2.77e-8}&\scriptsize{81.05}	 \vspace{0mm} \\
 \hline
 \multicolumn{8}{|c|}{\scriptsize{$\mathcal{T}=\mathbb{R}^{20\times20\times20\times20\times20}$, $r=(2,2,2,2,2)$}} \\
 \hline
 \multicolumn{4}{|c|}{\scriptsize{$sr=0.3$}} & \multicolumn{4}{c|}{\scriptsize{$sr=0.6$}} \\
 \hline
 \scriptsize{SALM-LRTC}	&\scriptsize{72}&\scriptsize{9.46e-9}&\scriptsize{246.45}& \scriptsize{SALM-LRTC} &\scriptsize{35}&\scriptsize{5.88e-9}&\scriptsize{117.54}	 \vspace{-1.5mm} \\	
 \scriptsize{FP-LRTC}	&\scriptsize{520}&\scriptsize{5.25e-8}&\scriptsize{387.94}& \scriptsize{FP-LRTC} &\scriptsize{105}&\scriptsize{2.14e-8}&\scriptsize{84.64}		 \vspace{-1.5mm} \\
 \scriptsize{TENSOR-HC}	&\scriptsize{57}&\scriptsize{2.15e-8}&\scriptsize{272.88}& \scriptsize{TENSOR-HC} &\scriptsize{40}&\scriptsize{2.74e-8}&\scriptsize{187.79}	 \vspace{-1.5mm} \\	
 \scriptsize{ADM-CON}	&\scriptsize{175}&\scriptsize{1.60e-8}&\scriptsize{813.11}& \scriptsize{ADM-CON} &\scriptsize{76}&\scriptsize{1.08e-8}&\scriptsize{375.13}		 \vspace{-1.5mm} \\
 \scriptsize{ADM-TR(E)}	&\scriptsize{381}&\scriptsize{1.97e-7}&\scriptsize{1194.78}& \scriptsize{ADM-TR(E)} &\scriptsize{192}&\scriptsize{2.19e-8}&\scriptsize{656.17}	 \vspace{0mm} \\
 \hline
 \end{tabular}
 \end{table}

 Different problem settings are used to test the algorithms. The order of tensors varies from three to five, and
 we also vary the multiliear-rank and the sampling ratio $sr$. Table \textbf{1} reports the results for random low
 multilinear-rank tensor completion problems without noise. In the table, we report the average iterations, the average
 relative error, and the average time (in seconds) of 10 runs. As can be seen from Table \textbf{1}, SALM-LRTC outperforms
 other algorithms in most cases. Although FP-LRTC converges fastest for the problems with high $sr$ and low $n$-rank
 (e.g., $\mathcal{T}=\mathbb{R}^{100\times100\times50}$, $r=(10,10,5)$, $sr=0.6$), SALM-LRTC also results nearly the same
 computation time and relative errors as those of FP-LRTC. For other problems, SALM-LRTC is robust and
 converges faster than others obviously. Especially, it takes SALM-LRTC no more than 40 iterations on the average to
 solve the problems with $sr=0.6$ in our experiments. In addition, the solutions for almost all the problems in our experiments
 computed by SALM-LRTC are more accurate than those delivered by other algorithms. Specifically, the relative errors of SALM-LRTC
 are smaller than $1 \times 10^{-8}$ in most cases.

 \begin{figure}[H]
 \centering
 \subfigure[]{\includegraphics[height=5.8cm]{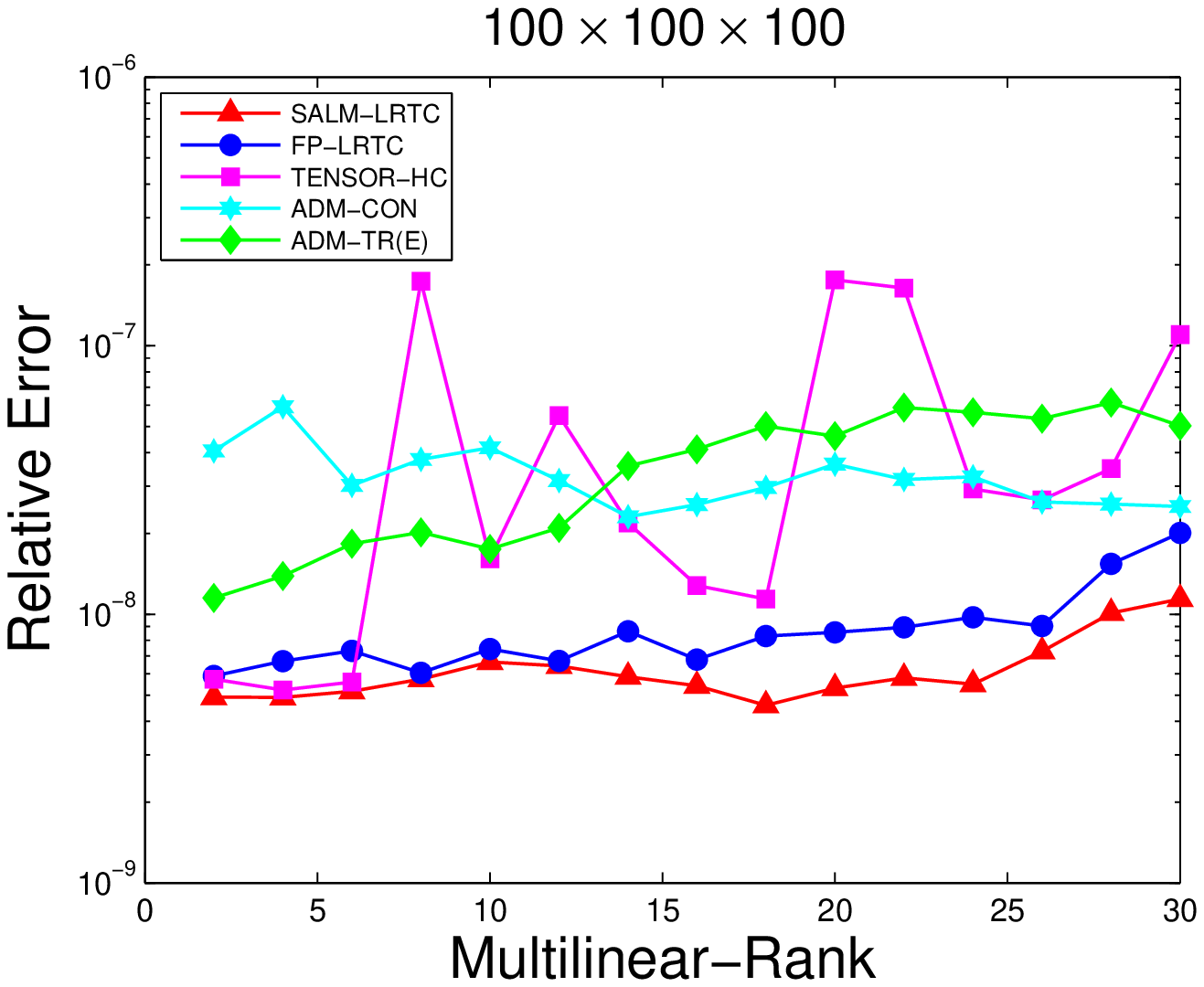}}
 \subfigure[]{\includegraphics[height=5.8cm]{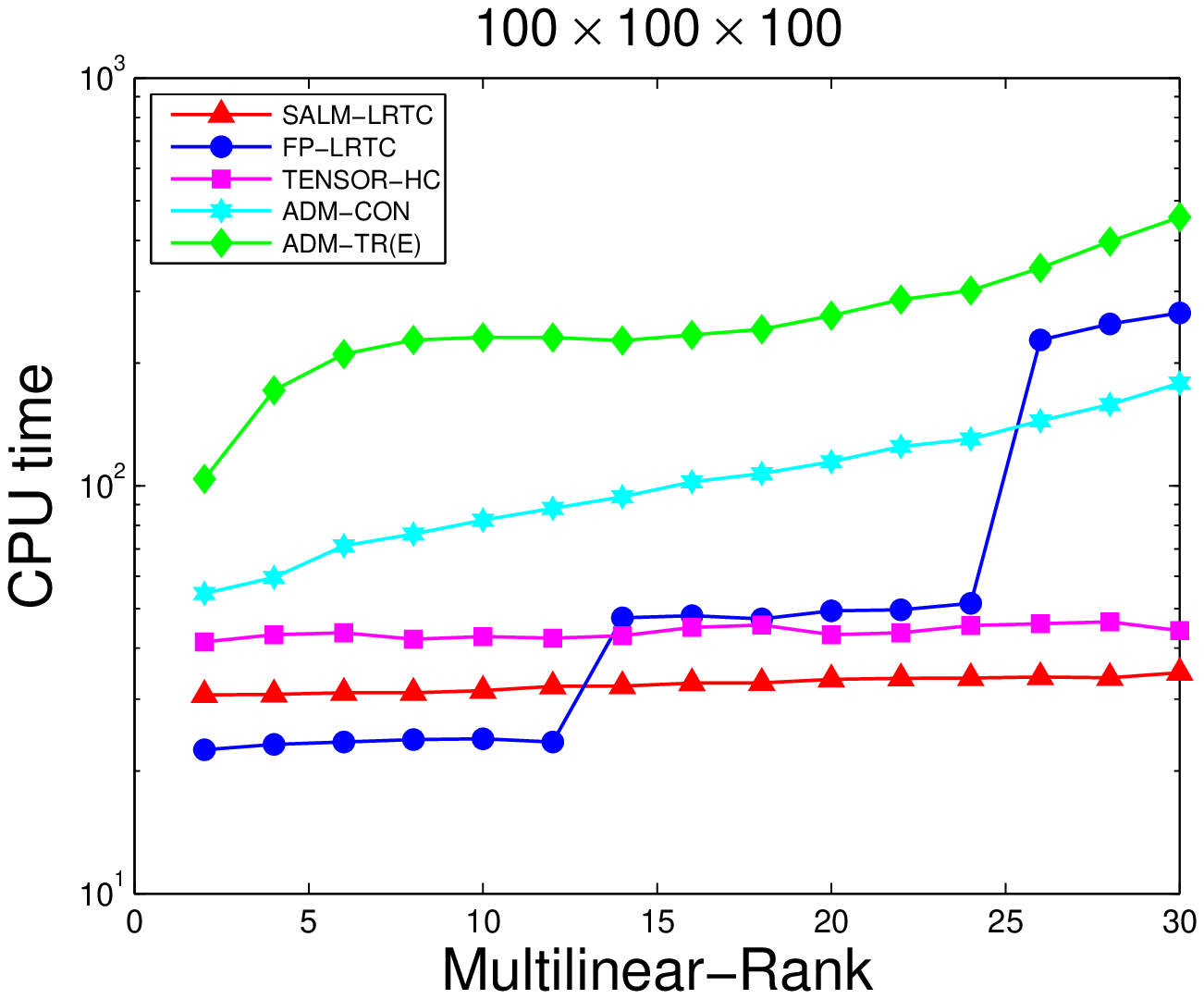}}
 \caption{Recovery results of $100 \times 100 \times 100$ tensors with $sr=0.5$ and different multilinear-ranks by SALM-LRTC, FP-LRTC, TENSOR-HC, ADM-CON and ADM-TR(E). (a) relative error; (b) CPU time in seconds. }
 \end{figure}

 In Figure \textbf{1}, we fix the tensor size ($100 \times 100 \times 100$) and compare the first four different
 algorithms (ADM-TR(E) is poorer than other algorithms obviously by Table \textbf{1}) for created noiseless problems
 with different multilinear-ranks $(r, r, r)$ (here we set $r_1 =r_2 =r_3=r$ for convenience). The sampling ratio is set to $sr=0.5$. Figure \textbf{1} shows that all the algorithms are effective for this set of problems. We can also see that SALM-LRTC costs relatively constant time and is faster than others for most problems. Moreover, SALM-LRTC is more robust and always has more accurate solution than that
 of other algorithms.

 We further study SALM-LRTC to solve random low multilinear-rank tensor completion problems with noisy data.
 Table \textbf{2} presents the numerical performance. In the table, we report the mean of NRMSEs, iterations and execution times over 10 independent trials. Note that a different $\sigma$ gives a different noise level. Then, we set $\sigma=0.02$ and $\sigma=0.04$, respectively. From the results, we can see that the NRMSEs of all the algorithms
 are smaller than the noise level in the given data. Especially, SALM-LRTC and TENSOR-HC have nearly the same
 performance, which are better than others.

 \begin{table}[H]
 \centering \tabcolsep 6pt
 \begin{tabular}{|llcr|llcr|}
 \multicolumn{8}{c}{\scriptsize{\textbf{Table 2.} Numerical results on low multilinear-rank tensor completion with noise}} \vspace{1mm} \\
 \hline
 \scriptsize{Method} & \scriptsize{iter} & \scriptsize{NRMSE} & \scriptsize{time} & \scriptsize{Method} & \scriptsize{iter} & \scriptsize{NRMSE} & \scriptsize{time} \\
 \hline
 \multicolumn{8}{|c|}{\scriptsize{$\mathcal{T}=\mathbb{R}^{50\times50\times50}$, $r=(9,9,3)$, $sr=0.3$}}  \\
 \hline
 \multicolumn{4}{|c|}{\scriptsize{$\sigma=0.02$}} & \multicolumn{4}{c|}{\scriptsize{$\sigma=0.04$}} \\
 \hline
 \scriptsize{SALM-LRTC}	&\scriptsize{40}&\scriptsize{1.24e-2}&\scriptsize{2.32}& \scriptsize{SALM-LRTC} &\scriptsize{40}&\scriptsize{2.11e-2}&\scriptsize{2.35}	 \vspace{-1.5mm}\\	
 \scriptsize{FP-LRTC}	&\scriptsize{500}&\scriptsize{1.29e-2}&\scriptsize{11.08}& \scriptsize{FP-LRTC} &\scriptsize{500}&\scriptsize{2.13e-2}&\scriptsize{11.25}		 \vspace{-1.5mm}\\
 \scriptsize{TENSOR-HC}	&\scriptsize{34}&\scriptsize{9.25e-3}&\scriptsize{2.97}& \scriptsize{TENSOR-HC} &\scriptsize{28}&\scriptsize{1.74e-2}&\scriptsize{2.48}	 \vspace{-1.5mm}\\	
 \scriptsize{ADM-CON}	&\scriptsize{106}&\scriptsize{1.11e-2}&\scriptsize{8.57}& \scriptsize{ADM-CON} &\scriptsize{140}&\scriptsize{2.13e-2}&\scriptsize{11.23}		 \vspace{-1.5mm}\\
 \scriptsize{ADM-TR(E)}	&\scriptsize{240}&\scriptsize{1.35e-2}&\scriptsize{27.99}& \scriptsize{ADM-TR(E)} &\scriptsize{314}&\scriptsize{2.15e-2}&\scriptsize{34.65}	 \vspace{0mm}\\
 \hline
 \multicolumn{8}{|c|}{\scriptsize{$\mathcal{T}=\mathbb{R}^{50\times50\times50}$, $r=(9,9,3)$, $sr=0.6$}}   \\
 \hline
 \multicolumn{4}{|c|}{\scriptsize{$\sigma=0.02$}} & \multicolumn{4}{c|}{\scriptsize{$\sigma=0.04$}} \\
 \hline
 \scriptsize{SALM-LRTC}	&\scriptsize{37}&\scriptsize{6.80e-3}&\scriptsize{2.04}& \scriptsize{SALM-LRTC} &\scriptsize{31}&\scriptsize{1.26e-2}&\scriptsize{1.69}	 \vspace{-1.5mm}\\	
 \scriptsize{FP-LRTC}	&\scriptsize{105}&\scriptsize{7.48e-3}&\scriptsize{2.44}& \scriptsize{FP-LRTC} &\scriptsize{105}&\scriptsize{1.37e-2}&\scriptsize{2.37}		 \vspace{-1.5mm} \\
 \scriptsize{TENSOR-HC}	&\scriptsize{26}&\scriptsize{9.84e-3}&\scriptsize{2.30}& \scriptsize{TENSOR-HC} &\scriptsize{20}&\scriptsize{1.89e-2}&\scriptsize{1.72}	 \vspace{-1.5mm} \\	
 \scriptsize{ADM-CON}	&\scriptsize{38}&\scriptsize{9.32e-3}&\scriptsize{3.13}& \scriptsize{ADM-CON} &\scriptsize{69}&\scriptsize{1.76e-2}&\scriptsize{5.68}		 \vspace{-1.5mm} \\
 \scriptsize{ADM-TR(E)}	&\scriptsize{445}&\scriptsize{6.77e-3}&\scriptsize{58.18}& \scriptsize{ADM-TR(E)} &\scriptsize{534}&\scriptsize{1.29e-2}&\scriptsize{76.14}	 \vspace{0mm} \\
 \hline
 \multicolumn{8}{|c|}{\scriptsize{$\mathcal{T}=\mathbb{R}^{20\times20\times20\times20\times20}$,
 $r=(2,2,2,2,2)$, $sr=0.3$}} \\
 \hline
 \multicolumn{4}{|c|}{\scriptsize{$\sigma=0.02$}} & \multicolumn{4}{c|}{\scriptsize{$\sigma=0.04$}} \\
 \hline
 \scriptsize{SALM-LRTC}	&\scriptsize{27}&\scriptsize{7.75e-3}&\scriptsize{86.85}& \scriptsize{SALM-LRTC} &\scriptsize{25}&\scriptsize{1.11e-2}&\scriptsize{84.95}	 \vspace{-1.5mm} \\	
 \scriptsize{FP-LRTC}	&\scriptsize{500}&\scriptsize{8.61e-3}&\scriptsize{436.99}& \scriptsize{FP-LRTC} &\scriptsize{500}&\scriptsize{1.01e-2}&\scriptsize{440.56}		 \vspace{-1.5mm} \\
 \scriptsize{TENSOR-HC}	&\scriptsize{18}&\scriptsize{8.58e-3}&\scriptsize{80.14}& \scriptsize{TENSOR-HC} &\scriptsize{20}&\scriptsize{1.59e-2}&\scriptsize{92.28}	 \vspace{-1.5mm} \\	
 \scriptsize{ADM-CON}	&\scriptsize{287}&\scriptsize{8.19e-3}&\scriptsize{1435.63}& \scriptsize{ADM-CON} &\scriptsize{359}&\scriptsize{1.26e-2}&\scriptsize{1802.40}		 \vspace{-1.5mm} \\
 \scriptsize{ADM-TR(E)}	&\scriptsize{389}&\scriptsize{5.84e-3}&\scriptsize{1162.46}& \scriptsize{ADM-TR(E)} &\scriptsize{689}&\scriptsize{1.25e-2}&\scriptsize{2015.95}	 \vspace{0mm} \\
 \hline
 \multicolumn{8}{|c|}{\scriptsize{$\mathcal{T}=\mathbb{R}^{20\times20\times20\times20\times20}$, $r=(2,2,2,2,2)$, $sr=0.6$}} \\
 \hline
 \multicolumn{4}{|c|}{\scriptsize{$\sigma=0.02$}} & \multicolumn{4}{c|}{\scriptsize{$\sigma=0.04$}} \\
 \hline
 \scriptsize{SALM-LRTC}	&\scriptsize{24}&\scriptsize{5.27e-3}&\scriptsize{82.61}& \scriptsize{SALM-LRTC} &\scriptsize{24}&\scriptsize{8.79e-3}&\scriptsize{83.23}	 \vspace{-1.5mm} \\	
 \scriptsize{FP-LRTC}	&\scriptsize{105}&\scriptsize{5.43e-3}&\scriptsize{95.15}& \scriptsize{FP-LRTC} &\scriptsize{105}&\scriptsize{9.03e-3}&\scriptsize{97.22}		 \vspace{-1.5mm} \\
 \scriptsize{TENSOR-HC}	&\scriptsize{15}&\scriptsize{9.04e-3}&\scriptsize{70.87}& \scriptsize{TENSOR-HC} &\scriptsize{20}&\scriptsize{1.74e-2}&\scriptsize{99.48}	 \vspace{-1.5mm} \\	
 \scriptsize{ADM-CON}	&\scriptsize{227}&\scriptsize{8.65e-3}&\scriptsize{1145.63}& \scriptsize{ADM-CON} &\scriptsize{336}&\scriptsize{1.55e-2}&\scriptsize{1700.12}		 \vspace{-1.5mm} \\
 \scriptsize{ADM-TR(E)}	&\scriptsize{334}&\scriptsize{5.15e-3}&\scriptsize{1056.33}& \scriptsize{ADM-TR(E)} &\scriptsize{455}&\scriptsize{7.54e-3}&\scriptsize{1402.68}	 \vspace{0mm} \\
 \hline
 \end{tabular}
 \end{table}

 \subsection{Image Simulation}

 In this part, we test the performance of SALM-LRTC on image inpainting \cite{bmgv2000,ng2006}. In fact, a color image can
 be represented as a third-order tensor. Then, if the image is of low multilinear-rank, or numerical low multilinear-rank,
 we can solve the image inpainting problem as a low multilinear-rank tensor completion problem. In other words,
 we can recover an image by solving problem \reff{problem5} with the sampling set $\Omega$ indexing non-missing pixels.

 In our test, for each image we remove entries in all the channels simultaneously (first two rows in Figure 2, $90.80\%$
 and $62.86\%$ known entries, respectively), or consider the case where entries are missing at random (last row in Figure 2,
 $30\%$ known entries). Figure 2 reports the original pictures, the input data tensor and the outcome of our algorithm.
 The relative errors of recovered results are 2.29e-2, 2.73e-2 and 1.09e-1, respectively.

 \begin{figure}[H]
 \centering
 \subfigure[]{\includegraphics[height=3.5cm,width=3.5cm]{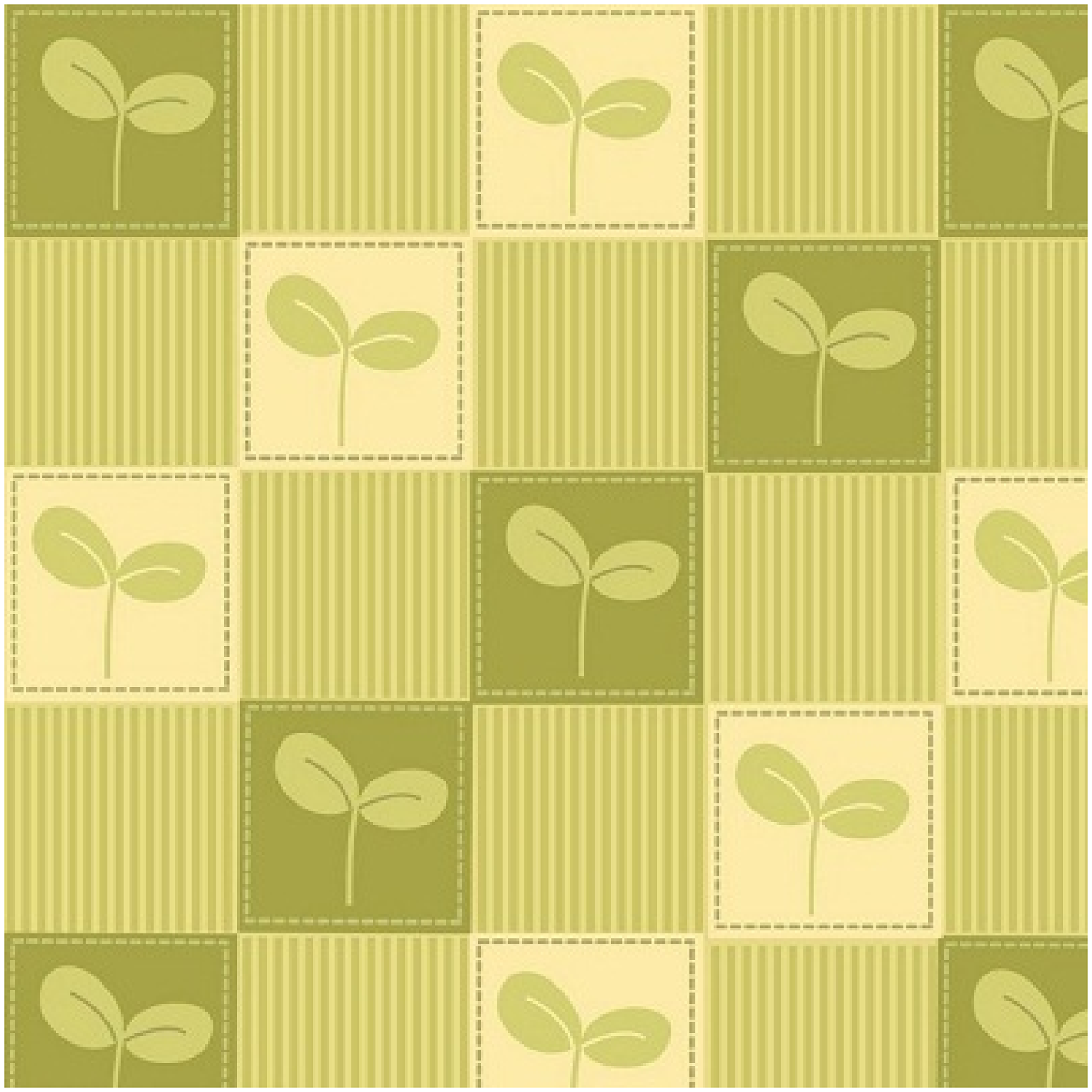}~
              \includegraphics[height=3.5cm,width=3.5cm]{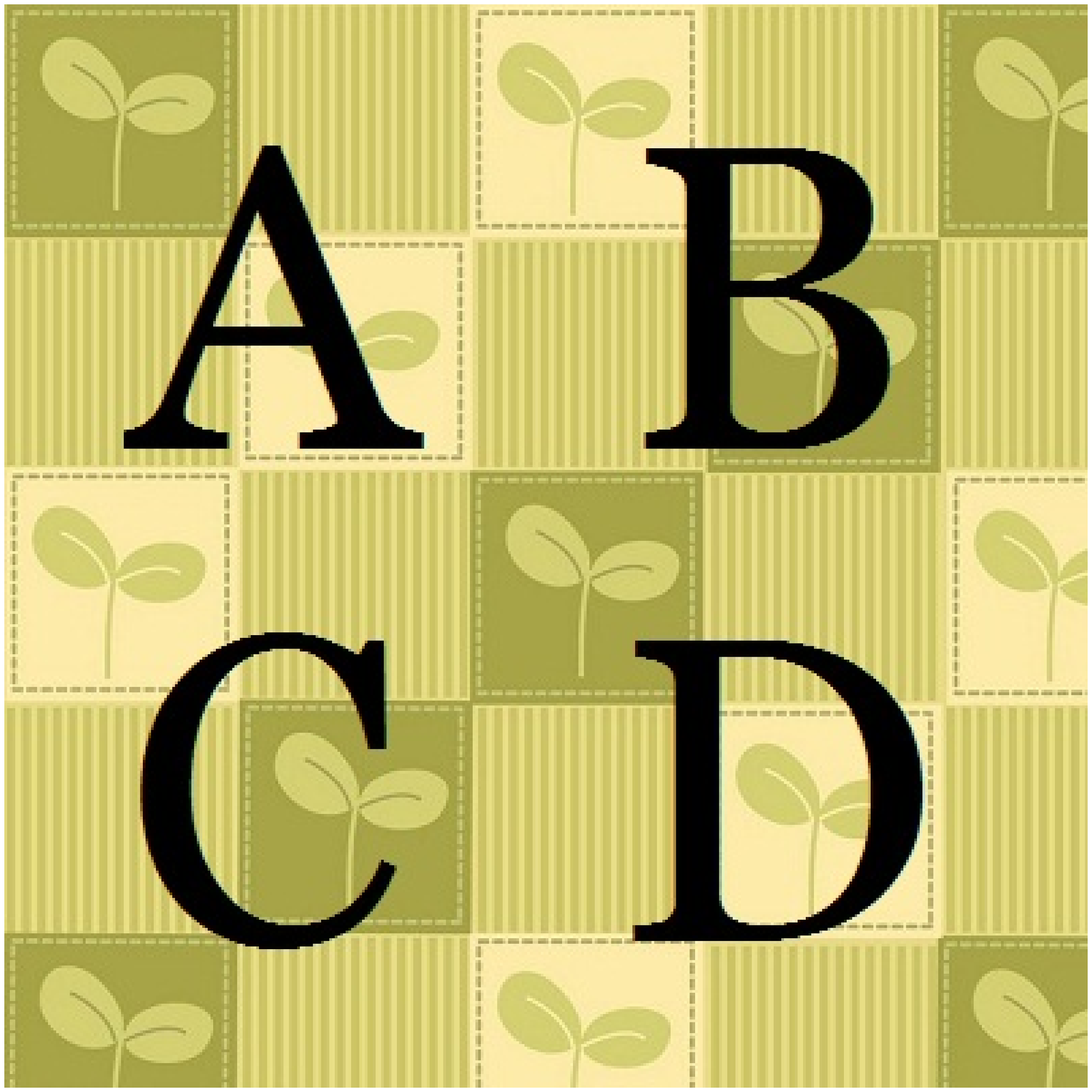}~
              \includegraphics[height=3.5cm,width=3.5cm]{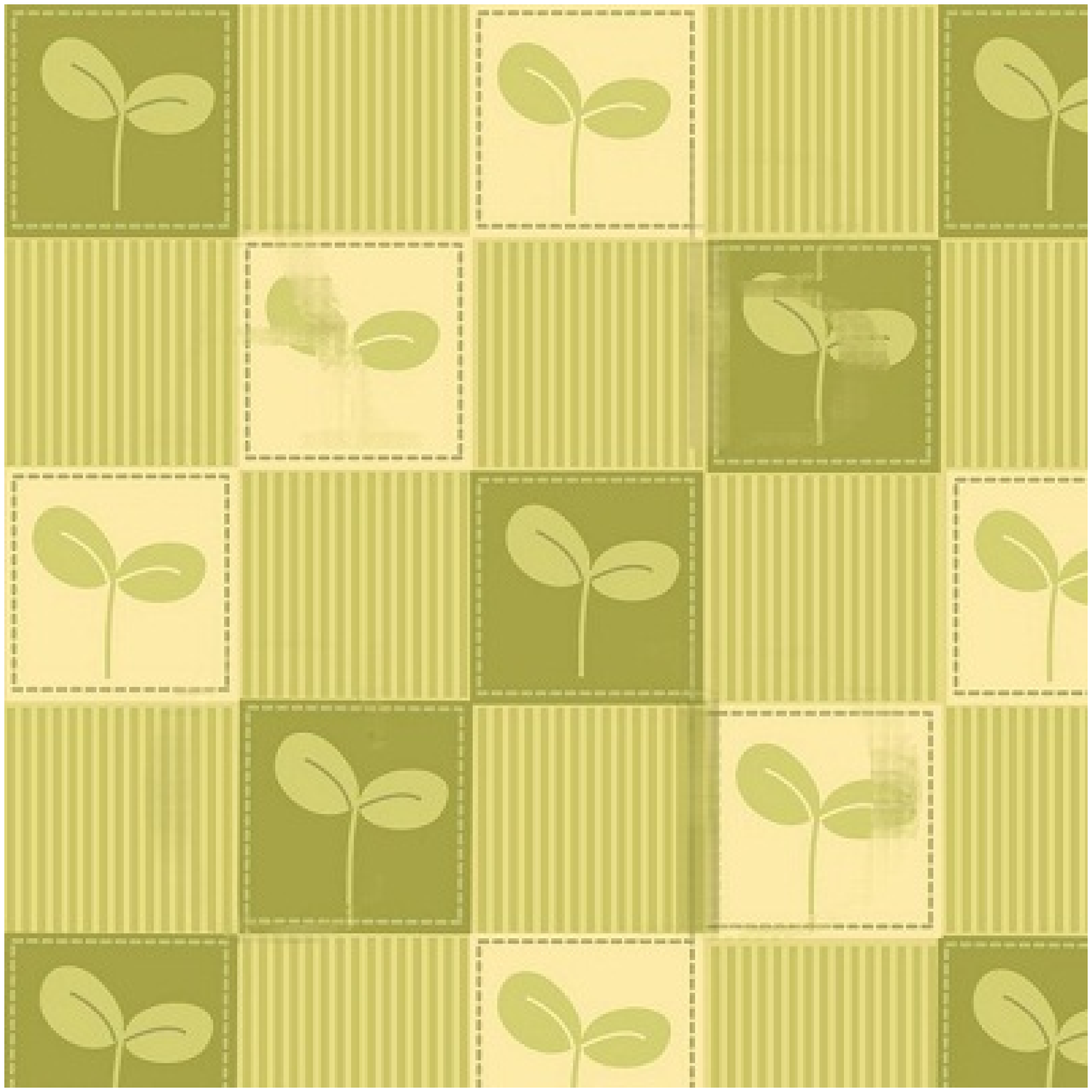}}
 \subfigure[]{\includegraphics[height=3.5cm,width=3.5cm]{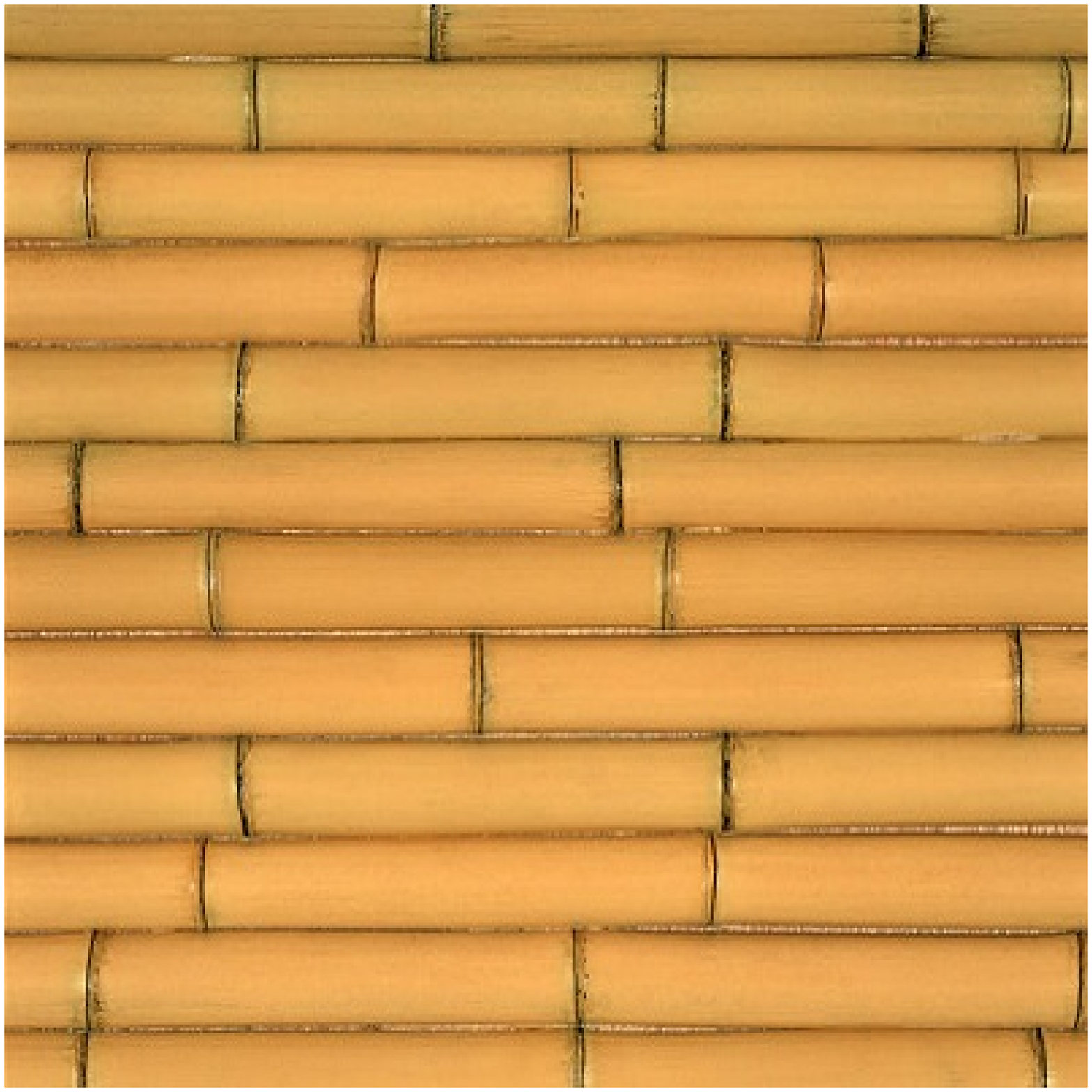}~
              \includegraphics[height=3.5cm,width=3.5cm]{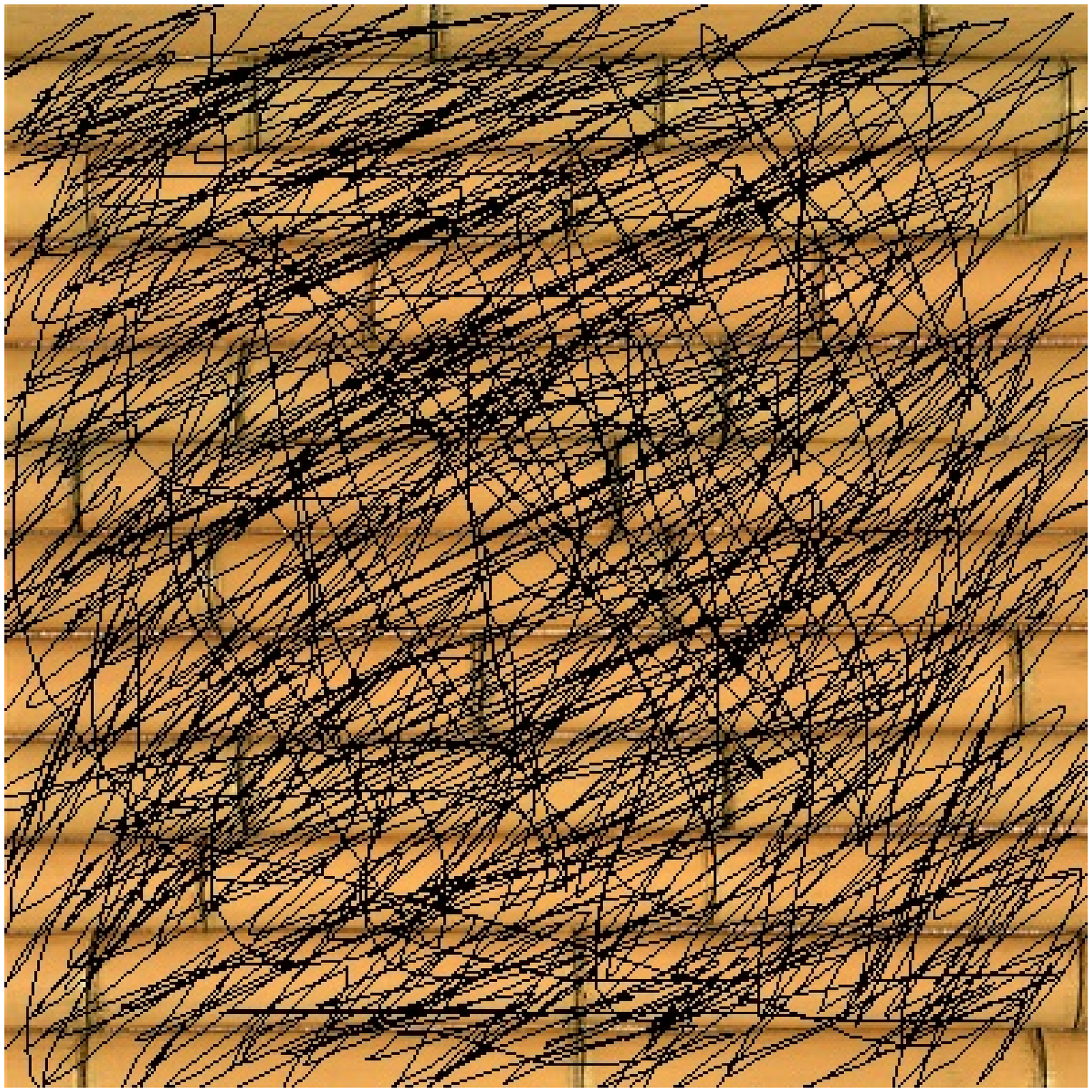}~
              \includegraphics[height=3.5cm,width=3.5cm]{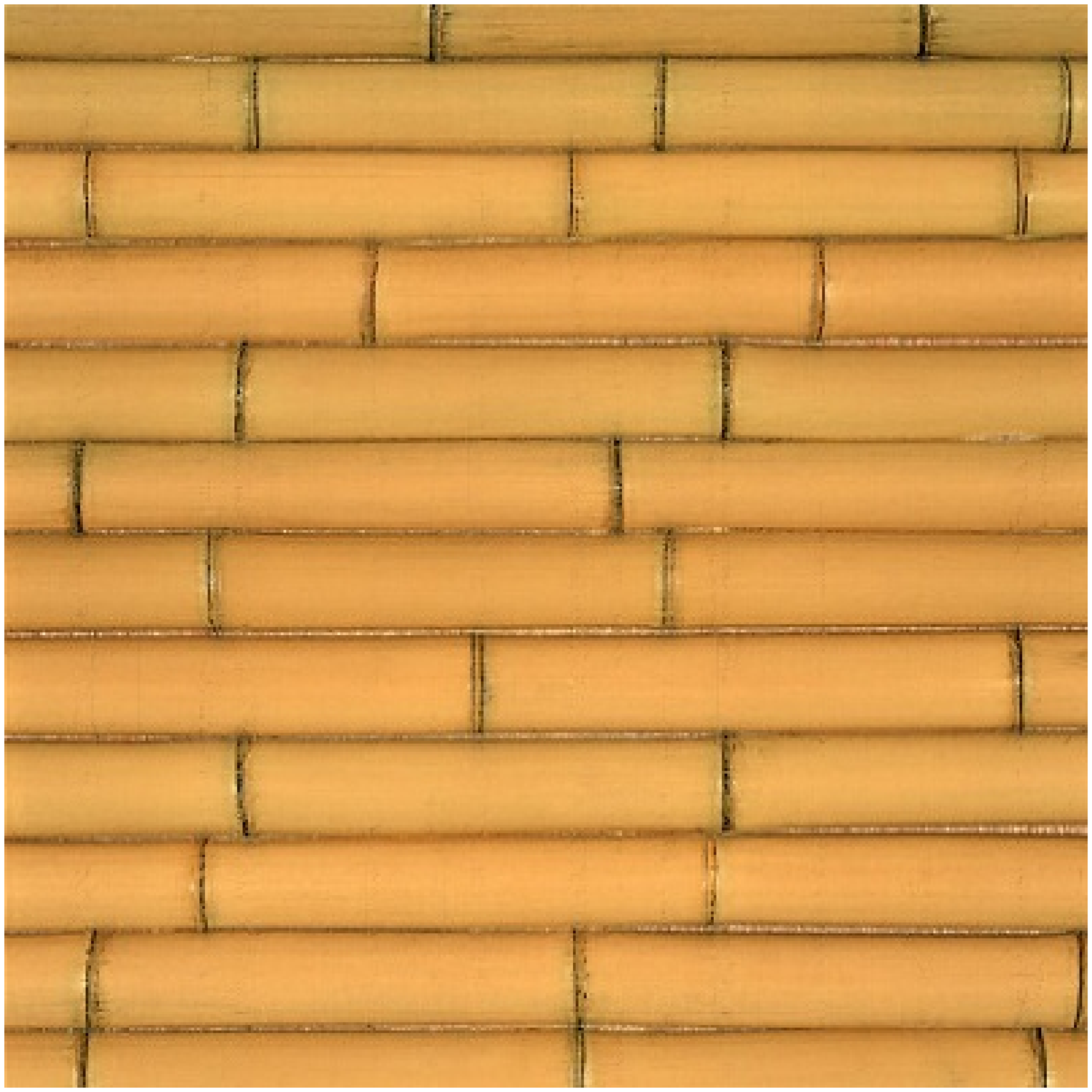}}
 \subfigure[]{\includegraphics[height=3.5cm,width=3.5cm]{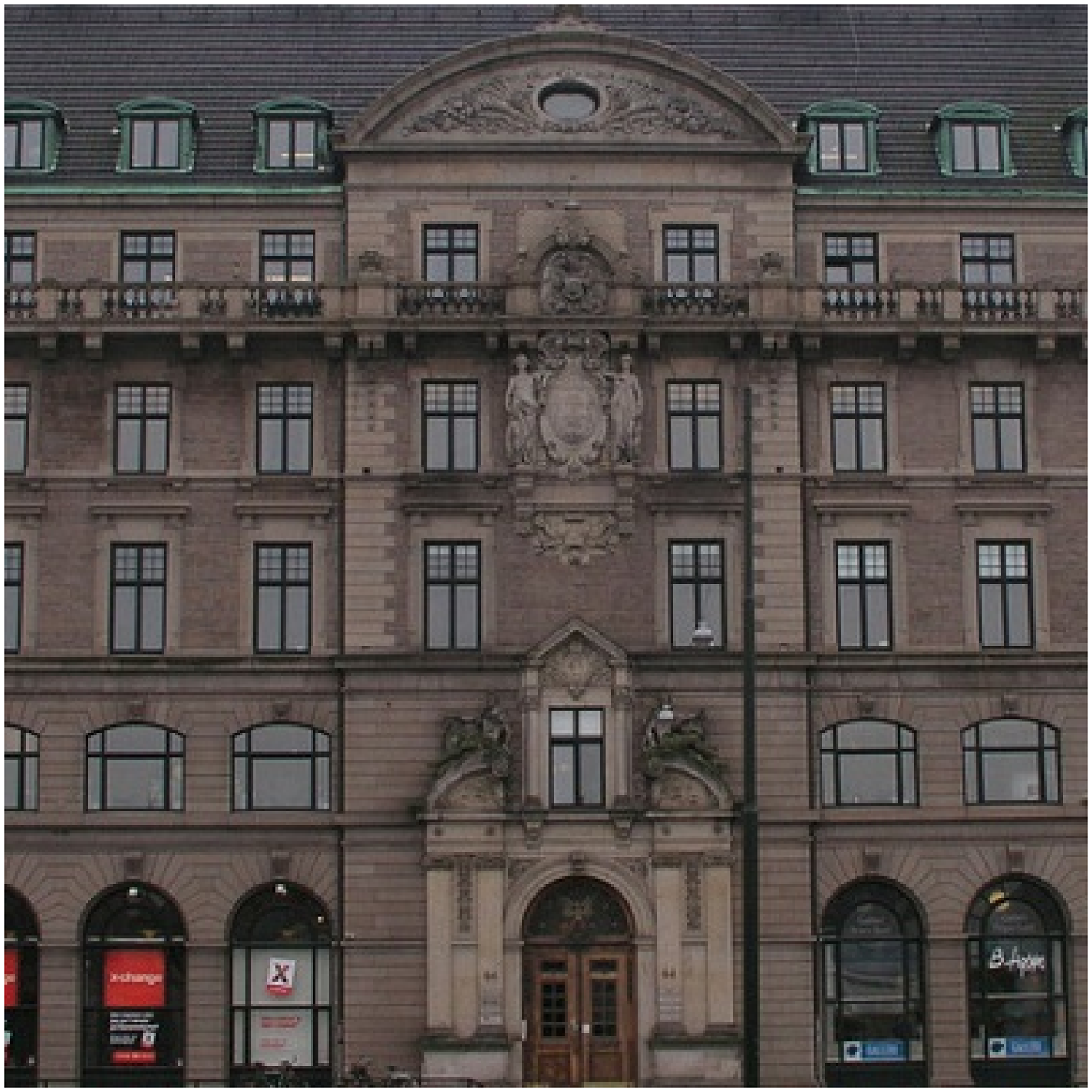}~
              \includegraphics[height=3.5cm,width=3.5cm]{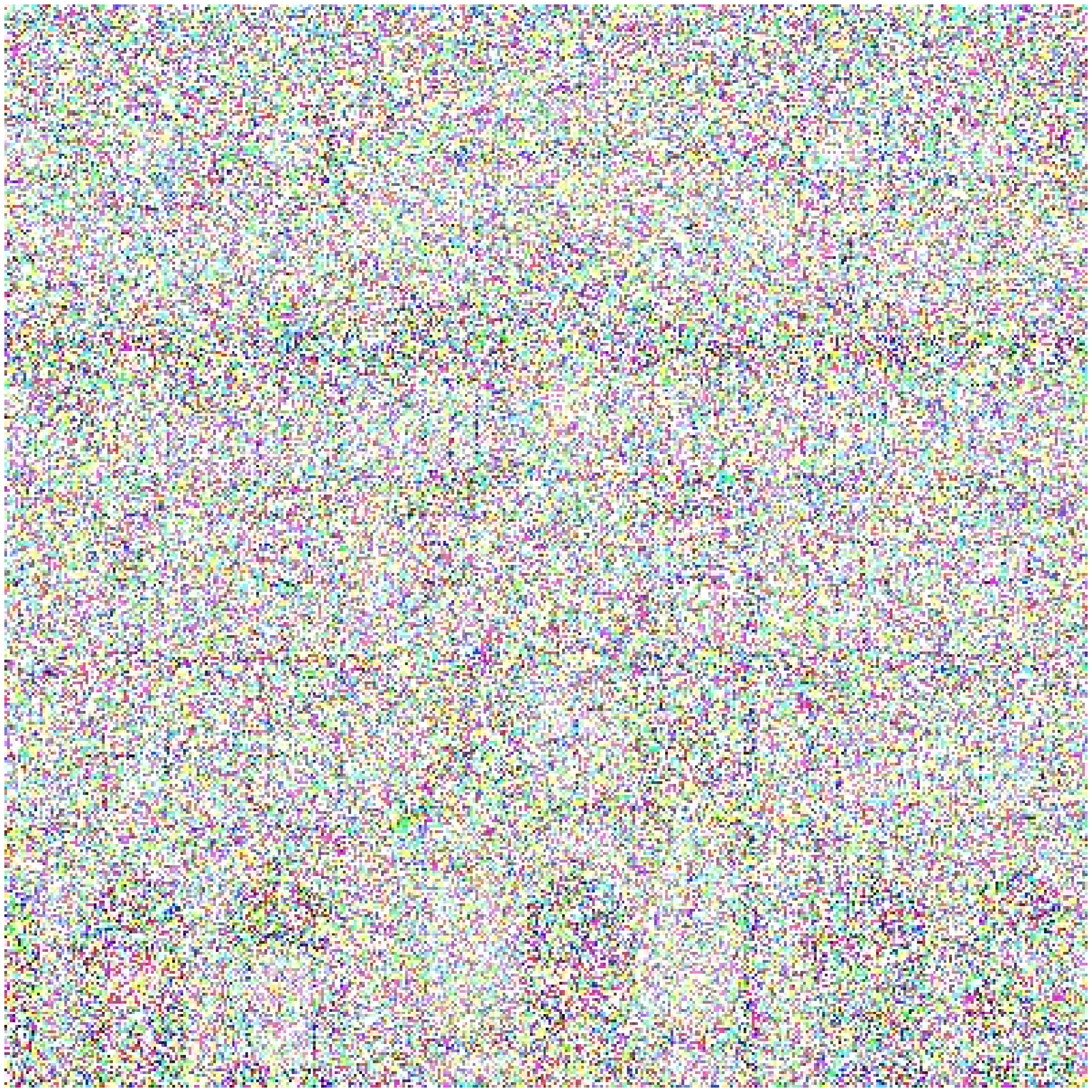}~
              \includegraphics[height=3.5cm,width=3.5cm]{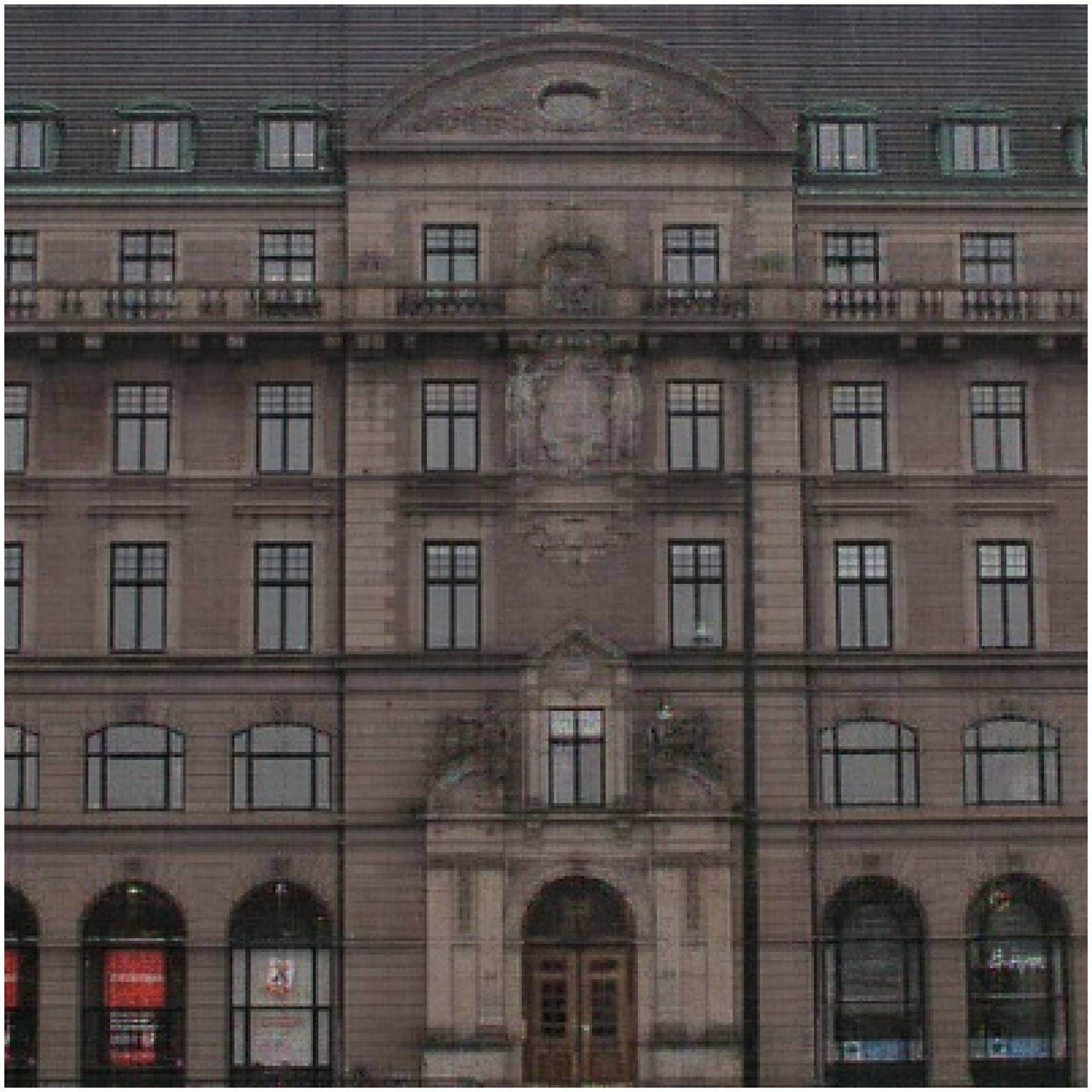}}
 \caption{Original images (left column); Input to the algorithm (middle column); The recovered result by SALM-LRTC (right column).}
 \end{figure}

 \section{ Conclusions }
 In this paper we focused on low multilinear-rank tensor recovery problem, and adopted variable splitting
 technique and convex relaxation technique to transform it into a tractable constrained optimization problem,
 which can be solved by the classical ALM directly. Taking advantage of well structure emerging in the
 transformed model, an easily implemented algorithm is developed based on the classical ALM.
 Some preliminary numerical results show the efficiency and robustness of the proposed algorithm.

 It is interesting to investigate how to make the best use of tensor structure to improve the relative model and
 algorithm. We believe that the augmented-Lagrangian-type framework can yield more robust and effective methods for
 more general tensor optimization problems based on the full utilization of the desired structure. Moreover,
 the nonconvex sparse optimization problems and the related algorithms in vector or matrix space have been widely
 discussed in the literature \cite{zxfh2012,mzy2013,ln2011,yzy2013}. It is worth investigating the nonconvex model
 in the tensor space.
 \vspace{6mm}

 \noindent {\bf Acknowledgments}\vspace{2mm}

 \noindent We would like to thank Silvia Gandy for sending us the code of ADM-TR(E), and thank
 Marco Signoretto for sending us the code of TENSOR-HC.



\begin{thebibliography}{1}

 \bibitem{tb2008}
 T.G. Kolda and B.W. Bader, ``Tensor decompositions and applications,'' {\em SIAM Rev.}, vol. 51,
 pp. 457-464, 2009.
 \vspace{-2mm}

 \bibitem{bmgv2000}
 M. Bertalm\'{\i}o, G. Sapiro, V. Caselles and C. Ballester, ``Image inpainting,'' {\em Proceedings
 of SIGGRAPH 2000}, New Orleans, USA, 2000.
 \vspace{-2mm}

 \bibitem{ng2006}
 N. Komodakis and G. Tziritas, ``Image completion using global optimization." CVPR, pp. 417-424, 2006.
 \vspace{-2mm}

 \bibitem{kr2007}
 T. Korah, C. Rasmussen, ``Spatiotemporal inpainting for recovering texture maps of occluded building facades." {\em IEEE Trans. Image Process.}, vol. 16(9), pp. 2262-2271, 2007.
 \vspace{-2mm}

 \bibitem{mhs2001}
 M. Fazel, H. Hindi and S. Boyd, ``A rank minimization heuristic with application to minimum order
 system approximation,'' {\em Proceedings of the American Control Conference (Arlington, VA, June
 2001)}, vol. 6, 4734-4739, 2001.
 \vspace{-2mm}

 \bibitem{bmp2010}
 B. Recht, M. Fazel and P. A. Parrilo, ``Guaranteed minimum-rank solutions of linear matrix equations
 via nuclear norm minimization,'' {\em SIAM Rev.}, vol. 52, pp. 471-501, 2010.
 \vspace{-2mm}

 \bibitem{zl2009}
 Z. Liu and L. Vandenberghe, ``Interior-point method for nuclear norm approximation with application to system identification."
 {\em SIAM J. Matrix Anal. Appl.}, vol. 31(3), pp. 1235-1256, 2009.
 \vspace{-2mm}

 \bibitem{ney1995}
 N. Linial, E. London, and Y. Rabinovich, ``The geometry of graphs and some of its algorithmic applications". {\em Combinatorica},
 vol. 15, pp. 215-245, 1995.
 \vspace{-2mm}

 \bibitem{j1990}
 J.H{\aa}stad, ``Tensor rank is NP-complete,'' {\em J. Algorithms}, vol. 11, pp. 644-654, 1990.
 \vspace{-2mm}

 \bibitem{eb2009}
 E.J. Cand\`{e}s and B. Recht, ``Exact matrix completion via convex optimization,'' {\em Found. Comput.
 Math.}, vol. 9, pp. 717-772, 2009.
 \vspace{-2mm}

 \bibitem{bwb2011}
 B. Recht, W. Xu and B. Hassibi, ``Null space conditions and threshlods for rank minimization," {\em Math. Program.},
 vol. 127, pp. 175-202, 2011.
 \vspace{-2mm}

 \bibitem{b2011}
 B. Recht, ``A simpler approach to matrix completion,'' {\em J. Mach. Learn. Res.}, vol. 12, pp.
 3413-3430, 2011.
 \vspace{-2mm}

 \bibitem{ydk2012}
 Y.J. Liu, D.F. Sun and K.C. Toh, ``An implementable proximal point algorithmic framework for nuclear norm minimization,"
 {\em Math. Program.}, vol. 133(1-2), pp. 399-436, 2012.
 \vspace{-2mm}

 \bibitem{bmp2007}
 M. Fazel, ``Matrix Rank Minimization with Applications,'' PhD thesis, Stanford University, 2002.
 \vspace{-2mm}

 \bibitem{sdl2009}
 S.Q. Ma, D. Goldfarb and L.F. Chen, ``Fixed point and Bregman iterative methods for matrix rank
 minimization,'' {\em Math. Program.}, vol. 128, pp. 321-353, 2011.
 \vspace{-2mm}

 \bibitem{jez2008}
 J.F. Cai, E.J. Cand\`{e}s and Z.W. Shen, ``A singular value thresholding algorithm for matrix
 completion,'' {\em SIAM J. Optim.}, vol. 20, No.4, pp. 1956-1982, 2010.
 \vspace{-2mm}

 \bibitem{ks2009}
 K.C. Toh and S.W. Yun, ``An accelerated proximal gradient algorithm for nuclear norm regularized
 linear least squares problems,'' {\em Pac. J. Optim.}, vol. 6, pp. 615-640, 2010.
 \vspace{-2mm}

 \bibitem{jppj2009}
 J. Liu, P. Musialski, P. Wonka and J.P. Ye, ``Tensor completion for estimating missing values in
 visual data,'' in {\em IEEE Int. Conf. Computer Vision (ICCV)}, Kyoto, Japan, pp. 2114-2121, 2009.
 \vspace{-2mm}

 \bibitem{sbi2011}
 S. Gandy, B. Recht and I. Yamada, ``Tensor completion and low-$n$-rank tensor recovery via convex
 optimization,'' {\em Inv. Probl.}, vol. 27, 025010 (19pp), 2011.
 \vspace{-2mm}

 \bibitem{mqlj2011}
 M. Signoretto, Q. Tran Dinh, L. De Lathauwer, J.A.K. Suykens, ``Learning with tensors: a framework
 based on convex optimization and spectral regularization," to appear in {\em Mach. Learn.}, 2013.
 \vspace{-2mm}

 \bibitem{yhs2012}
 L. Yang, Z.H. Huang and X.J. Shi, ``A fixed point iterative method for low $n$-rank tensor pursuit,"
 {\em IEEE Trans. Signal Process.}, vol. 61(11), pp. 2952-2962, 2013.
 \vspace{-2mm}

 \bibitem{mlj2010}
 M. Signoretto, L. De. Lathauwer and J.A.K. Suykens, ``Nuclear norms for tensors and their use for
 convex multilinear estimation,'' {\em Internal report 10-186, ESAT-SISTA, K.U. Leuven, Leuven,
 Belgium}, Lirias number: 270741, 2010.
 \vspace{-2mm}

 \bibitem{mrbj2011}
 M. Signoretto, R. Plas, B. Moor, and J. Suykens, ``Tensor versus matrix completion: A comparison with application to spectral data,"
 {\em IEEE Signal Process. Lett.}, vol. 18(7), pp. 403-406, 2011.
 \vspace{-2mm}

 \bibitem{rkh2011}
 R. Tomioka, K. Hayashi and H. Kashima, ``Estimation of low-rank tensors via convex optimization,''
 Arxiv preprint arXiv:1010.0789v2, 2011.
 \vspace{-2mm}

 \bibitem{rtkh2011}
 R. Tomioka, T. Suzuki, K. Hayashi, and H. Kashima, ``Statistical performance of convex tensor
 decomposition," Advances in Neural Information Processing Systems (NIPS) 24. 2011, Granada, Spain.
 \vspace{-2mm}

 \bibitem{mz2012}
 M. Zhang, Z.H. Huang, ``Exact recovery conditions for the low-$n$-rank tensor recovery problem via
 its convex relaxation," submitted to {\em J. Optim. Theory Appl.}, Revised 2013.
 \vspace{-2mm}

 \bibitem{mknr2011}
 M.E. Kilmer, K. Braman, N. Hao, R.C. Hoover, `` Third order tensors as
 operators on matrices: a theoretical and computational framework with applications in
 imaging," {\em SIAM J. Matrix Anal. Appl.}, vol. 34(1), pp. 148-172, 2013
 \vspace{-2mm}

 \bibitem{lzsj2013}
 L. Yang, Z.H. Huang, S.L. Hu and J.Y. Han, ``An iterative algorithm for third-order tensor multi-rank minimization,"
 submitted to {\em Comput. Optim. Appl.}, Revised 2013. \vspace{-7mm}

 \bibitem{ty2011}
 M. Tao and X. M. Yuan, ``Recovering low-rank and sparse components of matrices from incomplete
 and noisy observations," {\em SIAM J. Optim.}, vol. 21, pp. 57-81, 2011.
 \vspace{-2mm}

 \bibitem{hty2012}
 B.S. He, M. Tao and X. M. Yuan, ``Alternating direction method with Gaussian back substitution for
 separable convex programming," {\em SIAM J. Optim.}, vol. 22(2), pp. 313-340, 2012.
 \vspace{-2mm}

 \bibitem{lcm2009}
 Z. Lin, M. Chen and Y. Ma, ``The augmented lagrange multiplier method
 for exact recovery of corrupted low-rank matrices," {\em UIUC Technical Report
 UILU-ENG-09-2215}, 2009.
 \vspace{-2mm}

 \bibitem{m1969}
 M. Hestenes, ``Multiplier and gradient methods,"{\em J. Optim. Theory Appl.}, vol. 4, pp. 303-320, 1969.
 \vspace{-2mm}

 \bibitem{mj1969}
 M.J.D. Powell, A method for nonlinear constraints in minimization problems, in Optimization,
 R. Fletcher, ed., Academic Press, New York, 1969, pp. 283¨C298.
 \vspace{-2mm}

 \bibitem{r1970}
 R. T. Rockafellar, Convex Analysis. Princeton University Press, Princeton, 1970.
 \vspace{-2mm}

 \bibitem{b1992}
 M. W. Berry, ``Large-scale sparse singular value decompositions,'' {\em The International Journal of
 Supercomputer Applications}, vol. 6, pp. 13-49, 1992.
 \vspace{-2mm}

 \bibitem{propack}
 R. M. Larsen. PROPACK-software for large and sparse svd calculations available at http://sun.stanford.
 edu/srmunk/PROPACK/.
 \vspace{-2mm}

 \bibitem{zxfh2012}
 Z.B. Xu, X.Y Chang, F.M. Xu and H. Zhang, ``$L_{1/2}$ regularization: A thresholding representation theory and a fast solver,"
 {\em IEEE Trans. Neural Networ. Learn. Syst.}, vol. 23(7), pp. 1013-1027, 2012.
 \vspace{-2mm}

 \bibitem{mzy2013}
 M. Zhang, Z.H. Huang and Y. Zhang, ``Restricted $p$-isometry properties of nonconvex matrix recovery,"
 {\em IEEE Trans. Inform. Theory}, vol. 59(7), pp. 4316-4323, 2013.
 \vspace{-2mm}

 \bibitem{ln2011}
 L.C. Kong and N.H. Xiu, ``Exact low-rank matrix recovery via nonconvex schatten $p$-minimization,"
 {\em Asia. Pac. J. Oper. Res.}, vol. 30(3), 1340010(13pages), 2013.
 \vspace{-2mm}

 \bibitem{yzy2013}
 Y.F. Li, Y.J. Zhang and Z.H. Huang, ``Reweighted nuclear norm minimization algorithm for low rank matrix recovery,"
 submitted to {\em J. Comput. Appl. Math.}, Revised 2013.



 \end{thebibliography}
\end{document}